\documentclass[12pt]{amsart}

\usepackage{amssymb} 
\usepackage{amsfonts}
\usepackage{eucal}
\usepackage{amsmath} 
\usepackage{amscd} 
\usepackage{latexsym} 
\usepackage{graphicx}
\usepackage{multirow}

\def\F{\mathbb{F}}

\def\C{\mathbb{C}}
\def\Z{\mathbb{Z}}
\def\Q{\mathbb{Q}}

\def\Fp{\F_p}

\def\p{\mathbf{p}}
\def\q{\mathbf{q}}
\def\r{\mathbf{r}}
\def\s{\mathbf{s}}
\def\v{\mathbf{v}}
\def\w{\mathbf{w}}
\def\z{\mathbf{z}}
\def\e{\mathbf{e}}

\newcommand{\smin}{\text{-}}
\newcommand{\smins}{\text{--}}
\newcommand{\splus}{\text{+}}

\newcommand{\recray}[4]{ \mbox{\tiny $
\begin{array}{r@{\hspace{4pt}}r@{\hspace{4pt}}r@{\hspace{4pt}}r}
#1
\end{array}
\left[ 
\begin{array}{r@{\hspace{4pt}}r@{\hspace{4pt}}r@{\hspace{4pt}}r} 
#2 
\end{array}
 \left| 
\begin{array}{r@{\hspace{4pt}}r@{\hspace{4pt}}r@{\hspace{4pt}}r} 
#3
\end{array}
  \right|  
\begin{array}{r@{\hspace{4pt}}r@{\hspace{4pt}}r@{\hspace{4pt}}r} 
#4  
\end{array}
\right] $ }
 }

\newcommand{\wvec}[1]{W\hspace{-0.035 in}\left( #1\right)}
\newcommand{\myenspace}{}
\newcommand{\Spec}{\operatorname{Spec}}

\newcommand{\mathcolon}{\mathrel{\mathop:}}

\usepackage{xypic}
\usepackage{graphpap}
\usepackage{url}

\include{biblio}

\newtheorem{theorem}{Theorem}[section]
\newtheorem{corollary}[theorem]{Corollary}
\newtheorem{lemma}[theorem]{Lemma} 
\newtheorem{proposition}[theorem]{Proposition} 
\newtheorem*{theoremnonum}{Theorem}

\theoremstyle{definition}
\newtheorem{definition}[theorem]{Definition}

\theoremstyle{remark}
\newtheorem{remark}[theorem]{Remark}

\theoremstyle{remark}
\newtheorem{example}[theorem]{Example}

\title{Elliptic Nets and Elliptic Curves}
\date{\today} 
\author{Katherine E. Stange} 
\address{Department of Mathematics, Simon Fraser University,
8888 University Drive, Burnaby, BC, Canada V5A 1S6,
and Pacific Institute for the Mathematical Sciences,
200 1933 West Mall, Vancouver, BC, Canada V6T 1Z2}
\email{kestange@pims.math.ca}

\thanks{This work was supported by NSERC Awards PGS D2 331379 and PDF 373333.}

\keywords{elliptic curve, elliptic divisibility sequence, recurrence sequence}

\subjclass[1991/2000]{Primary 11G05, 11G07, 11B37, Secondary 11B39, 14H52}

\begin{document}

\begin{abstract}
An elliptic divisibility sequence is an integer recurrence sequence associated to an elliptic curve over the rationals together with a rational point on that curve.  In this paper we present a higher-dimensional analogue over arbitrary base fields.  Suppose $E$ is an elliptic curve over a field $K$, and $P_1, \ldots, P_n$ are points on $E$ defined over $K$.  To this information we associate an $n$-dimensional array of values in $K$ satisfying a nonlinear recurrence relation.  Arrays satisfying this relation are called \emph{elliptic nets}.  We demonstrate an explicit bijection between the set of elliptic nets and the set of elliptic curves with specified points.  We also obtain Laurentness/integrality results for elliptic nets.
\end{abstract}
\maketitle
\setcounter{tocdepth}{1}
\tableofcontents

\section*{Introduction}

An \emph{elliptic divisibility sequence} is an integer sequence $W_n$ satisfying %the recurrence relation
\begin{equation}
\label{eqn:ward}
W_{n+m}W_{n-m}= W_{n+1}W_{n-1}W_m^2-W_{m+1}W_{m-1}W_n^2.
\end{equation}
This definition was introduced by Morgan Ward in 1948 \cite{War}.  Let $\Psi_n(x,y)$ denote the $n$-th division polynomial associated to an elliptic curve (the $n$-th division polynomial vanishes at the $n$ torsion points).  Ward showed that division polynomials satisfy the recurrence \eqref{eqn:ward} and furthermore that all elliptic divisibility sequences have the form
\[
W_n = \lambda^{n^2-1}\Psi_n(x,y)
\]
for some constant $\lambda$, elliptic curve (or singular cubic) and point $P = (x,y)$ on the curve.  This rich structure has led to number theoretic results~\cite{Aya, EveMclWar, Ing2, Sil5, Sil4, Swa}; applications to Hilbert's 10th problem \cite{CorZah, EisEve, Poo2}; to integrable systems \cite{Hon}; and to
cryptography~\cite{ChuChu, Shi, Sta3}.  For a bibliography, see \cite[Chapter 10]{EvePooShpWar}.
\par
There have been several attempts to generalise this theory.  Van der Poorten and Swart study \emph{translated elliptic divisiblity sequences} \cite{Swa, Poo, PooSwa}.  Mazur and Tate generalise division polynomials to arbitrary endomorphisms in the $p$-adic setting \cite{MazTat}, and Streng uses their definition to generalise to the endomorphism ring of an elliptic curve with complex multiplication \cite{Str}.  Elliptic divisibility sequences are closely related to the denominators of the multiples $[n]P$ of a fixed point $P$; questions have been asked about the collection of denominators of the linear combinations $[n]P + [m]Q$ by Everest, Miller and Stephens \cite{EveMilSte}.  The hope of defining `higher rank' elliptic divisibility sequences via a recurrence relation was discussed in correspondence by Elkies, Propp and Somos \cite{Rob}.  
\par
The primary purpose of this paper is to generalise from integer sequences to multi-dimensional arrays with values in any field, which we call \emph{elliptic nets}.  A substantial part of the difficulty lies in finding the correct recurrence %relation %% POPS THE PAGE
and defining a generalised division polynomial.
\par
We define an \emph{elliptic net} to be a function $W: A \rightarrow R$ from a finite rank free abelian group $A$ to an integral domain $R$ satisfying the properties that $W(0)=0$ and that
\begin{multline}
\label{eqn:netintro}
\wvec{p + q + s} \wvec{p- q}\wvec{r+ s}\wvec{r} \\
+ \wvec{q+ r+ s}\wvec{ q- r}\wvec{ p+ s}\wvec{ p} \\
+ \wvec{ r + p+ s} \wvec{ r- p} \wvec{ q+ s}\wvec{ q} = 0
\end{multline}
for all $p,q,r,s \in A$.  If $A=R=\Z$, this is an equivalent definition of an elliptic divisibility sequence (this is not immediately obvious, but it is a consequence of results in this paper).  By the \emph{rank} of an elliptic net we shall mean the rank of $A$ (this bears no relation to the \emph{rank of apparition} defined by Ward for elliptic divisibility sequences \cite{War}).  Section \ref{sec: edsedn} covers the basic definitions and gives examples.
\par
Our primary interest is the relationship between elliptic curves and elliptic nets.
\begin{theoremnonum}[Main Theorem - Introductory Version]{}
For each field $K$ and integer $n$, there is an explicit bijection of sets
\begin{equation*}
\xymatrix{
\left\{
{
\begin{array}{l}
\mbox{scale equivalence classes} \\
\mbox{of non-degenerate elliptic} \\ \mbox{nets }
W: \Z^n \rightarrow K
\end{array}
}
\right\} \ar[d] \\
\left\{
{
\begin{array}{l}
\mbox{tuples }(C, P_1, \ldots , P_n )\mbox{ where $C$ is a cubic}\\
\mbox{curve in Weierstrass form defined over $K$,}\\
\mbox{considered modulo unihomothetic changes}\\
\mbox{of variables, and such that }\{P_i\} \in C_{ns}(K)^n \\
\mbox{is appropriate}
\end{array}
}
\right\} \ar[u]
}
\end{equation*}
\end{theoremnonum}

For a description of the relevant terminology, see Sections \ref{sec: arbitraryfields} (appropriate), \ref{sec:scalenormal} (scale equivalent, non-degenerate) and \ref{sec:homosing} (unihomothetic).  See Theorem \ref{thm:big} for a more detailed statement.  The isomorphism itself is described explicitly in Definition \ref{defn: WEP} (depending on Theorem \ref{thm: arbitraryfield}) and Theorem \ref{thm: curvefromnetn}.  For ranks 1 and 2, explicit formul{\ae} can be found in Propositions~\ref{prop: 2explicit}, \ref{prop: curvefromnet1} and \ref{prop: curvefromnet2}. For an example, see Figure \ref{figure: ellnet66}.
\par
The other main aspect of elliptic nets studied in this paper is Laurentness.  These results are needed for the proof of the main theorem, but are of independent interest.  One property of elliptic divisibility sequences of particular interest is that they are integer sequences:  if the sequence begins $1, a, b, ac, \ldots (a,b,c \in \Z)$, then it will consist entirely of integers \cite{War}.  This result has been studied in the more general framework of the `Laurent phenomenon' of Fomin and Zelevinsky \cite{FomZel}.  
\par
Laurentness results are found in Section \ref{sec: induction}, which is devoted to the inductive structure of elliptic nets:  how some terms are determined by others via the recurrence relation.  We define a universal ring $\mathcal{W}_A$ for elliptic nets on $A$, such that elliptic nets $W: A \rightarrow R$ are in bijection with homomorphisms $\mathcal{W}_A \rightarrow R$.  We obtain results on the structure of this ring, and in turn, these imply integrality results.  See Theorems~\ref{thm: 1laurent} ($n=1$), \ref{thm: 2laurent} ($n=2$) and \ref{thm: nlaurent} ($n \ge 3$).  The proofs in this section are elementary but somewhat tedious.

Sections \ref{sec: overc} and \ref{sec: netpoly} define the higher rank generalisation of division polynomials called \emph{net polynomials}:  rational functions on the $n$-fold product $E^n$ of an elliptic curve $E$, which vanish on tuples $(P_1, \ldots, P_n)$ satisfying a linear relation $[v_1]P_1 + \ldots + [v_n]P_n = \mathcal{O}$ for fixed coefficients $v_i$.  In Section~\ref{sec: overc}, we work with the complex uniformization of an elliptic curve defined over $\C$.  Section~\ref{sec: netpoly} generalises the definition to arbitrary fields by analysing the arithmetic properties of net polynomials.  The main result here is Theorem~\ref{thm: valuation}.  

The last three sections describe the bijection in the main theorem.  Section~\ref{sec: arbitraryfields} makes explicit the production of an elliptic net from any cubic Weierstrass curve using the net polynomials.  Section~\ref{sec: curvesfromnets} determines exactly those cubic curves which produce a given elliptic net.  Finally, Section~\ref{sec: curvenet} puts together the results of the previous sections to prove the main theorem, stated in its full form as Theorem \ref{thm:big}.

{\bf Computer software.}  The explicit isomorphism described in this paper has been implemented for Pari/GP in ranks 1 and 2.  Scripts are available at \url{http://math.katestange.net}.

{\bf Acknowledgements.} I would like to thank my thesis advisor, Joseph Silverman, for many patient hours.  I would also like to thank Rafe Jones, Alf van der Poorten, and Jonathan Wise.   %This work was supported by NSERC Awards PGS-D2 331379-2006 and PDF 373333-2009.

%\pagebreak

\section{Elliptic nets}
\label{sec: edsedn}

The following definition is the subject of the paper.

\begin{definition}
\label{def: ellnet}
Let $A$ be a free finitely-generated abelian group, and $R$ be an integral domain.  An \emph{elliptic net} is any map $W: A \rightarrow R$ with
\begin{equation}
\label{eqn: ellrec0}
W(0)=0,
\end{equation}
and such that for all $p$, $q$, $r$, $s \in A$,
\begin{multline}
\label{eqn: ellrec}
\wvec{p + q + s} \wvec{p- q}\wvec{r+ s}\wvec{r} \\
+ \wvec{q+ r+ s}\wvec{ q- r}\wvec{ p+ s}\wvec{ p} \\
+ \wvec{ r + p+ s} \wvec{ r- p} \wvec{ q+ s}\wvec{ q} = 0.
\end{multline}
\end{definition}

%Taking $p=q=r=s=\mathbf{0}$ in \eqref{eqn: ellrec}, we have $3W(\mathbf{0})^4 = 0$.  Thus, if $R$ is not of characteristic three, then the condition that $W(\mathbf{0}) = 0$ is entailed by \eqref{eqn: ellrec}.  In fact, t
Functions $W: A \rightarrow R$ which satisfy \eqref{eqn: ellrec} but not \eqref{eqn: ellrec0} can only appear in characteristic $3$ (to see this, take $p=q=r=s=0$ in \eqref{eqn: ellrec}).  Any constant function in characteristic 3 is an example.  By definition, these are not elliptic nets.

We refer to the rank of $A$ as the \emph{rank} of the elliptic net.  Suppose that $B \subset A$ is a subgroup of $A$.  Then the restriction to $B$ of an elliptic net $W: A \rightarrow R$ is also an elliptic net.  We refer to this elliptic net as \emph{the subnet associated to $B$} and write $W |_B : B \rightarrow R$.

\begin{example}\label{ex: nets} Let $R$ be an integral domain.  The following are elliptic nets.
\begin{enumerate}
\item The \emph{zero net} $W: \Z^n \rightarrow R$ defined by $W(\v)=0$ for all $\v$.
\item The identity map $W_{id}: \Z \rightarrow \Z$ given by $W(v)=v$.
\item Let $W': \mathbb{Z} \rightarrow R$ be an elliptic net.  Then for each $1 \leq i \leq n$, we may define $W_i: \Z^n \rightarrow R$ by $W_i(v_1,\ldots,v_n) =  W'(v_i)$, and this will also be an elliptic net.
\item More generally, if $W: A \rightarrow R$ is an elliptic net and $f: B \rightarrow A$ is a homomorphism of finitely generated free abelian groups, then $W \circ f : B \rightarrow R$ is also an elliptic net.
\item If $W: A \rightarrow R$ is an elliptic net and $g: R \rightarrow S$ is a homomorphism of integral domains, then $g \circ W : A \rightarrow S$ is also an elliptic net.
\item $W_{Leg}: \Z \rightarrow \Z$ given by $W(v)=\left(\frac{v}{3}\right)$, the Legendre symbol of $v$ over $3$.  This can be verified by a finite examination of cases; observe that at least one of~$p$, $q$, $r$, $p-q$, $q-r$, and $r-p$ is divisible by $3$.  See also \cite[p. 31]{War}.
\item $W_{Fib}: \Z \rightarrow \Z$ given by 
\[
W(v) = \left\{ \begin{array}{ll}
F_{2v} & v > 0 \\
-F_{2v} & v < 0 \\
0 & v = 0
\end{array} \right. .
\]
where $F_{2v}$ is the $2v$-th Fibonacci number.  One may verify this example using the closed form for terms of the Fibonacci sequence.  See also \cite[p. 31]{War}.
\item Figure \ref{figure: ellnet66} shows a portion of an elliptic net of rank 2 displayed as an array.  The origin is located at the term `0'.  This elliptic net arises from a certain curve and two points as described in Section \ref{sec: arbitraryfields}, Example \ref{example:net}.  Each axis forms an elliptic divisibility sequence, e.g. $0, 1, 1, -3, 11, 38, 249, \ldots$.
\end{enumerate}
\end{example}

\begin{figure}
\begin{tiny}
\caption{\small{Elliptic net associated to $y^2 + y = x^3 + x^2 - 2x$, $P= (0,0)$, $Q=(1,0)$ over $\Q$ (origin is at `0')}}
\label{figure: ellnet66}
\begin{equation*}
\!\!\!\!\!\!\!\!
\begin{array}{c}
 \\
 \\
 \\
 \\
 \\
 \\
% \\
% \\
% \\
% \\
 \uparrow \\
 Q
\end{array}\!\!\!\!
\begin{array}{l}
{\begin{array}{lllllllll}
%\hline
3269 & -2869 & {4335} & 5959 & {12016} & {-55287} & 23921 & 1587077 & {-7159461} \\
%\hline
-127 & -299 & {94} & 479 & 919 & -2591 & 13751 & {68428} & {424345} \\
%\hline
 -44 & -27 & -31 & 53 & {-33} & {-350} & 493 & {6627} & 48191 \\
%\hline
 -1 & -7 & {-5} & {8} & -19 & -41 & -151 & 989 & {-1466} \\
%\hline
 3 & -2 & 1 & {3} & -1 & -13 & {-36} & 181 & {-1535} \\
%\hline
1 & -1& 1 & 1 & {2} & {-5} & 7 & 89 & -149 \\
%\hline
 -1 & -1& {0} & 1 & 1 & {-3} & 11 & {38} & {249} \\
%\hline
-2 & -1 & -1 & 1 & -1 & -4 & 1 & 47 & 185 \\
%\hline
1 & -3 & -1 & 2 & -3 & -5 & -17& 63 & -184 \\
%\hline
%19 & -8 & 5 & 7 & 1 & -19 & -26 & 151 & -841\\
%\hline
%33 & -53 & 31 & 27 & 44 & -67 & 129 & 709 & 55 \\
%\hline
%-919 & -479 & -94 & 299 & 127 & -535 & 1187 & 3376 & 19061 \\
%\hline
%-12016 & -5959 & -4335 & 2869 & -3269 & -6522 & -3079 & 47987 & 87039 \\
%\hline
\end{array}} \\
P \rightarrow
\end{array}
\end{equation*}
\end{tiny}
\end{figure}

\section{Laurentness and integrality}
\label{sec: induction}

In this section we ask which terms of an elliptic net determine the others via the recurrence relation.  In the case of $n=1$, Ward showed that the terms $W(1), \ldots, W(4)$ sufficed to determine the rest of the net (unless too many of these terms were zero) \cite{War}.  Our method also demonstrates Laurentness and integrality results.  The main theorems of this section are used in Section \ref{sec: curvesfromnets}.

\subsection{Laurentness}

Let $I$ be a group, in additive notation, called the \emph{indexing group}, whose elements are called \emph{indices}.  To each $i \in I$, we associate the symbol $T_i$.  In what follows, the indexing group will be $I \cong \Z^n$ for some $n$.  

Consider the ideal $\mathcal{M}$ in the ring $\Z[T_i]_{i \in I}$ generated by $T_0$ and all polynomials
\begin{equation}
\label{eqn:recrel}
T_{p + q + s} T_{p-q}T_{r+s}T_{r}
+ T_{q+ r+s}T_{q-r}T_{p+s}T_{p}
+ T_{ r+p+ s} T_{r-p} T_{q+s}T_{q}
\end{equation}
(of the form \eqref{eqn: ellrec}) as $p, q, r, s$ range over $I$.  Polynomials of the form \eqref{eqn:recrel} will be called \emph{recurrence relations}.  
Consider the ring $\mathcal{W}_I$ obtained from $\Z[T_i]_{i \in I}/\mathcal{M}$ as a quotient by its own nilradical.  For each integral domain $R$, there is a bijection between elliptic nets $W: I \rightarrow R$ and homomorphisms $\mathcal{W}_I \rightarrow R$ (defined by taking $T_i \mapsto W(i)$).

Taking $p=q=i, r=s=0$ shows that $T_i^3(T_i+T_{-i}) \in \mathcal{M}$ for each $i \in I$.  In particular, $T_{-i}^3(T_i+T_{-i}) \in \mathcal{M}$ also.  Therefore, any prime ideal containing $\mathcal{M}$ contains $T_i + T_{-i}$; for if it did not, then it must contain $T_i$ and $T_{-i}$, a contradiction.  Therefore $T_{-i} = -T_i$ in $\mathcal{W}_I$.  This implies the following.

\begin{proposition}
\label{prop: netsymmzero}
Let $W: A \rightarrow R$ be an elliptic net.  Then $W(-z) = - W(z)$ for all $z \in A$.
\end{proposition}

The purpose of this section is to find a finite subset $0 \notin J \subset I$ such that the localisation $\mathcal{W}_I[T_i^{-1}]_{i \in J}$ is finitely generated as a $\Z$-algebra, and to give the generators.  (The localisation is not the trivial ring ($1=0$) by the existence of a homomorphism from it to $\Q$ given by Example \ref{ex: nets}, where one uses part (3) with $W'=W_{id}$ of part (2).)  From this we show that every $T_i$ can be expressed as a Laurent polynomial in integer coefficients in a finite number of terms $T_j$.  This implies that any elliptic net which does not take zero values at the $T_j$ is entirely determined by those values.  

To illustrate, consider the rank one case, which is essentially a result of Morgan Ward.

\begin{theorem}[{Ward, \cite[Theorem 4.1]{War}}]
\label{thm: 1laurent}
The ring $\mathcal{W}_\Z[T_1^{-1},T_2^{-1}]$ is generated as a $\Z$-algebra by the six elements
\[
T_1, \quad T_1^{-1}, \quad T_2, \quad T_2^{-1}, \quad T_3, \quad T_4.
\]
Furthermore, each $T_i$ is expressible as a $\Z$-coefficient polynomial in
\[
T_1, \quad T_1^{-1}, \quad T_2, \quad T_3, \quad T_4T_2^{-1}.
\]
%In other words, every element of $\mathcal{W}_\Z$ can be expressed as a polynomial with integer coefficients in these five variables.
In particular, let $W: \Z \rightarrow \Q$ be an elliptic net.  If $W(1)=1$, $W(2) \neq 0$, $W(i)$ is an integer for $i=2,3,4$, and $W(2)$ divides $W(4)$, then the elliptic net consists entirely of integers.
\end{theorem}

\begin{proof}
See \cite{War}, Theorem 4.1.  Recall that $T_{-n}=-T_n$, so it suffices to prove the first two statements for positive $n$.  Taking $(p,q,r,s) = (n+1,n,1,0)$ and $(n+1,n-1,1,0)$ respectively, in $\mathcal{W}_I$ we have
\begin{gather}
\label{eqn:2n1ind}
T_{2n+1}T_1^3 + T_{n-1}T_{n+1}^3 + T_{n+2}T_{-n}T_{n}^2 = 0,\\
\label{eqn:2ninduct}
T_{2n}T_2T_1^2 + T_nT_{n-2}T_{n+1}^2 + T_{n+2}T_{-n}T_{n-1}^2 = 0.
\end{gather}
The equations \eqref{eqn:2n1ind} and \eqref{eqn:2ninduct} prove the first statement by induction.  The base case consists of $0 \le n \le 4$; for $n > 4$, we have $2n > n+2$.

For even $i$, it can be shown by induction on \eqref{eqn:2ninduct} that $T_{i}$ is expressible as a $\Z$-coefficient polynomial in $T_1$, $T_1^{-1}$, $T_2$, $T_2^{-1}$, $T_3$, and $T_4$ in such a way that the combined degree of $T_2$ and $T_4$ in each monomial is positive.  For $i=2,4$ this is clear.
To complete the induction in general, observe that in \eqref{eqn:2ninduct}, each of the rightmost two terms is divisible by at least two $T_k$ where $k$ is even and $k < 2n$.

For even $i$, the second statement of the theorem concerning the expressibility of all $T_i$ in terms of $T_1$, $T_1^{-1}$, $T_2$, $T_3$ and $T_4T_2^{-1}$ follows from the observation of the previous paragraph.  The statement also holds for $i=1,3$.  Consequently, it holds for odd $i$ by induction on \eqref{eqn:2n1ind}.
\end{proof}

\subsection{Proofs by induction}
\label{sec:proofbyind}
The inductive proofs in this section will be based on the following definitions.  
Consider finite sets $S,~J~\subset~I$ where $0,i~\notin~S~\cup~J$.  We say that an index $i \in I$ is \emph{$S$-integrally implied by~$J$} if there exists a $\Z$-coefficient monomial $P(T_s)$ (in variables indexed by~$S$) and $\Z$-coefficient polynomial $Q(T_j)$ (in variables indexed by~$J$) such that
\begin{equation}
\label{eqn:imply}
T_iP(T_s)=Q(T_j)
\end{equation}
in $\mathcal{W}_I$.  A set $K \subset I$ is \emph{$S$-integrally implied by} the set $J$ if every index in $K$ is $S$-integrally implied by $J$.  

As an example (see Proposition \ref{prop: netsymmzero} and the paragraph which preceeds it), $-i$ is $S$-integrally implied by any $J$ containing $i$ (for any $S$).  In what follows, this fact will often be used tacitly.

A set $B \subset I$ \emph{is an $S$-integral baseset for $\mathcal{W}_I$} if all of $I$ is $S$-integrally implied by $B$.  If $B \subset I$ is an $S$-integral baseset, then each $T_i$ can be expressed as a polynomial with integer coefficients in the set of variables $\{T_b\}_{b \in B} \cup \{T_s^{-1}\}_{s \in S}$ (when considered in the appropriate localisation).

It is straightforward to verify that if $i$ is $S$-integrally implied by $J$ and every $j \in J$ is $S$-integrally implied by $J'$, then $i$ is $S$-integrally implied by $J'$.  To show that $B$ is an $S$-integral baseset for $I$, the proofs in this section show the following: for each index $i \in I$, there is a finite sequence $J_0 \subset J_1 \subset \cdots \subset J_n$ such that $B = J_0$, $i \in J_n$ and for each $1 \leq k \leq n$, $J_k$ is $S$-integrally implied by $J_{k-1}$.  At each stage, we show that each index of $J_i$ is $S$-integrally implied by $J_{i-1}$.  Recall that implication is simply the existence of an relation of the form \eqref{eqn:imply}, and in fact we simply give a relevant element of the form \eqref{eqn:recrel}.

These elements are cumbersome to write out.  For example, taking in the case $n=3$,
\[
\p = (1,0,0), \quad %\colvec{ 1 \\ 0 \\ 0 } ,\myenspace
\q = (0,1,0), \quad%\colvec{ 0 \\ 1 \\ 0 } ,\myenspace
\r = (0,0,1), \quad %\colvec{ 0 \\ 0 \\ 1 } ,\myenspace
\s = (0,0,0),%\colvec{ 0 \\ 0 \\ 0 } ,
\]
we obtain the element 
\begin{multline*}
T_{(1,1,0)}T_{(1,-1,0)}T_{(0,0,1)}T_{(0,0,1)}+\\
T_{(0,1,1)}T_{(0,1,-1)}T_{(1,0,0)}T_{(1,0,0)}+\\
T_{(1,0,1)}T_{(-1,0,1)}T_{(0,1,0)}T_{(0,1,0)}.
\end{multline*}
For this information, let us instead use a more convenient notation
\begin{equation}
\label{eqn: recarray}
\recray{
 1 & 0 & 0 & 0 \\
 0 & 1 & 0 & 0 \\
 0 & 0 & 1 & 0 
}{ 
 1 & 1 & 0 & 0 \\
 1 & \smin1 & 0 & 0 \\
 0 & 0 & 1 & 1 
 }{ 
 0 & 0 & 1 & 1 \\
 1 & 1 & 0 & 0 \\
 1 & \smin1 & 0 & 0 
 }{ 
 1 & \smin1 & 0 & 0 \\
 0 & 0 & 1 & 1 \\
 1 & 1 & 0 & 0 
  } \myenspace .
\end{equation}
In this notation, the columns to the left of the square braces correspond to the columns of $p$, $q$, $r$ and $s$, while the indices of the terms of the recurrence appear as the columns within the square braces.  

To demonstrate that an index $i$ is ($S$-integrally) implied by a set of indices $J$, it suffices to write down an appropriate such array.  Notice that any array of the form \eqref{eqn: recarray} is a recurrence if each row is a recurrence.  Therefore we may construct examples row-by-row.

The following definition will be useful for ordering inductions.

\begin{definition}
Let $$N(\v) = \max_{i=1,\ldots,n} |v_i|$$ be the \emph{sup-norm} of the vector $\v$.% or the term $W(\v)$.
\end{definition}

\subsection{Basesets for rank $2$}

For the rank two case, we require a lemma.

\begin{lemma}
\label{lemma: integrality}
The ring $\mathcal{W}_{\Z^2}[T_{(1,0)}^{-1}, T_{(0,1)}^{-1}, T_{(1,1)}^{-1}]$ is generated as a $\Z$-algebra by the elements
$$ \{ T_{\v} : N(\v) \leq 4 \} \cup \{ 
T_{(1,0)}^{-1}, T_{(0,1)}^{-1}, T_{(1,1)}^{-1} \}\myenspace . $$
\end{lemma}

\begin{proof}
Let $S = \{ (1,0), (0,1), (1,1) \}$ and $B = \{ \v \in \Z^2 : N(\v) \le 4 \}$.  This proof proceeds by induction on the sup-norm.  Trivially, any $\v$ with $N(\v) \le 4$ is $S$-integrally implied by $B$.  Let $N_0 > 4$ and suppose that all terms with indices with sup-norm less than $N_0$ are $S$-integrally implied by $B$.  Call the set of such indices $K_{N_0}$.  Suppose $\v$ is an index of sup-norm $N_0$.  We construct a recurrence demonstrating that $\v$ is $S$-integrally implied by $K_{N_0}$ row-by-row.  For $i = 1,2 $, define $w_i = \lceil\frac{v_i}{2}\rceil$.
%
%$$
%w_i = \left\{ \begin{array}{ll} v_i/2, & \myenspace v_i\text{ even} \\ (v_i + 1) / 2, & \myenspace v_i\text{ odd} \end{array} \right. \myenspace .
%$$%

{\bf Case I:  $\v$ has one odd entry and one even entry.}
For the odd entry, we use the row
\begin{equation*}
\recray{
w_i & w_i\smins 1 & 0 & 0
}{
v_i & 1 & 0 & 0
}{
w_i\smins 1 & w_i \smins 1 & w_i & w_i 
}{
w_i & \smins w_i & w_i\smins 1 & w_i\smins 1
}
\end{equation*}
For the even entry, we use the row
\begin{equation*}
\recray{
w_i & w_i & 1 & 0
}{
v_i & 0 & 1 & 1
}{
w_i\splus 1 & w_i \smins 1 & w_i & w_i 
}{
w_i\splus 1 & \smins w_i\splus 1 & w_i & w_i
}
\end{equation*}

{\bf Case II:  $\v$ has two odd entries.}
Use the rows
\begin{equation*}
\recray{
w_1 & w_1\smins 1 & 0 & 0 \\
w_2 & w_2 \smins 1 & 1 & 0 
}{
v_1 & 1 & 0 & 0\\
v_2 & 1 & 1 & 1
}{
w_1\smins 1 & w_1 \smins 1 & w_1 & w_1 \\
w_2 & w_2 \smins 2 & w_2 & w_2 
}{
w_1 & \smins w_1 & w_1\smins 1 & w_1\smins1 \\
w_2 \splus 1 & \smins w_2\splus 1 & w_2\smins1 & w_2\smins1
}
\end{equation*}

{\bf Case III:  $\v$ has two even entries.}
Use the rows
\begin{equation*}
\recray{
w_1 & w_1\smins 1 & 0 & 1 \\
w_2 & w_2 & 1 & 0 
}{
v_1 & 1 & 1 & 0\\
v_2 & 0 & 1 & 1
}{
w_1 & w_1 \smins 1 & w_1 \splus 1 & w_1 \\
w_2\splus1 & w_2 \smins 1 & w_2 & w_2
}{
w_1\splus1 & \smins w_1 & w_1 & w_1\smins1 \\
w_2 \splus 1 & \smins w_2\splus 1 & w_2 & w_2
}
\end{equation*}

For even $v_i$, either $|v_i| \leq 2$ or $|v_i| > 3$.  In the former case, $|w_i| + 1 \leq 2 < N_0$.  In the latter case, we have $|w_i| + 1 \leq ( |v_i| + 2 )/ 2 < |v_i| \leq N_0$.  For odd $\v_i$, either $|v_i| \leq 3$ or $|v_i| > 4$.  In the former case $|w_i| + 2 \leq 4 < N_0$.  In the latter case, we have $|w_i|+2 \leq ( |v_i| + 5 )/ 2 < |v_i| \leq N_0$.

Therefore all the vectors in the recurrence have sup-norm less than $N_0$ with the exception of $\v$.  In the monomial of $\v$ in the recurrence, the other indices are $(1,0)$, $(0,1)$ or $(1,1)$.  This demonstrates that $\v$ is $S$-integrally implied by $K_{N_0}$ and hence by $B$.
\end{proof}

\begin{theorem}
\label{thm: 2laurent}
The ring $\mathcal{W}_{\Z^2}[T_{(1,1)}^{-1},T_{(1,0)}^{-1},T_{(0,1)}^{-1}]$ is generated as a $\Z$-algebra by the eleven elements
\begin{gather*}
T_{(1,1)},\;\; T_{(1,0)},\;\; T_{(0,1)}, \;\; T_{(1,1)}^{-1},\;\; T_{(1,0)}^{-1},\;\; T_{(0,1)}^{-1},\\  T_{(2,1)},\;\;  T_{(1,2)}, \;\;  T_{(2,0)}, \;\; T_{(0,2)},\;\;  T_{(2,2)},
\end{gather*}
and the following identities hold:
\begin{align*}
T_{(1,-1)}T_{(1,1)}^3 &= 
T_{(1,0)}^3T_{(1,2)}-
T_{(0,1)}^3T_{(2,1)},\\
T_{(2,2)}T_{(1,-1)}T_{(1,0)}T_{(0,1)} &= T_{(1,1)}\left(T_{(0,2)}T_{(2,1)}T_{(1,0)} - T_{(0,1)}T_{(2,0)}T_{(1,2)} \right).
\end{align*}
\par
In particular, if $W: \Z^2 \rightarrow \Q$ is an elliptic net for which
\begin{enumerate}
\item $W(1,0) = W(0,1) = W(1,1) = 1$,
\item $W(2,0)$, $W(0,2)$, $W(1,2) \neq W(2,1)$ are integers, and
\item $W(1,2) - W(2,1)$ divides $W(0,2)W(2,1) - W(2,0)W(1,2)$,
\end{enumerate}
then all terms of the elliptic net are determined by these seven values and are integers.
\end{theorem}

\begin{proof}
The first and second stated identities are the recurrences
\begin{equation}
\label{eqn:1min1}
\recray{
 1 & 0 & 1 & 0 \\
 0 & 1 & 1 & 0 
 }{ 
  1 & 1 & 1 & 1 \\
 1 & \smin1 & 1 & 1 
  }{ 
  1 & \smin1 & 1 & 1 \\
 2 & 0 & 0 & 0 
  }{ 
  2 & 0 & 0 & 0 \\
 1 & 1 & 1 & 1 
   } \myenspace ,
\end{equation}
\begin{equation*}
\label{eqn: 2m2}
 \recray{
1 & 1 & \smin1 & 0 \\
 1 & 2 & 1 & \smin1 
}{ 
  2 & 0 & \smin1 & \smin1 \\
 2 & \smin1 & 0 & 1 
  }{ 
  0 & 2 & 1 & 1 \\
 2 & 1 & 0 & 1 
  }{ 
  0 & \smin2 & 1 & 1 \\
 1 & 0 & 1 & 2 
   } \myenspace .
\end{equation*}
%Notice that in \eqref{eqn: 2m2}, the calculation of $T_{(2,2)}$ requires division by $T_{(1,-1)}$.
%$$ T_{(1,-1)}T_{(1,1)}^3 = T_{(0,1)}^3T_{(2,1)} - T_{(1,0)}^3T_{(1,2)} \myenspace .$$

Let $S = \{ (1,0), (0,1), (1,1) \}$, and $B = \{ \v \in \Z^2 : N(\v) \leq 4 \}$.  By Lemma \ref{lemma: integrality}, it suffices to show that $B$ is $S$-integrally implied by the set
$$\{ (1,0), (0,1), (1,1), (2,0), (0,2), (2,1), (1,2), (2,2) \} \myenspace .$$
We list the relevant recurrences in order.  As each index is implied, it may be used to imply later indices.  It is assumed that as $(a,b)$ is implied, so is $(-a,-b)$.  To begin, the index $(1,-1)$ is implied by \eqref{eqn:1min1}.
\begin{equation*}
(2,-1):  
\recray{
 \smin1 & 0 & 1 & 1 \\
 1 & 1 & 0 & 0 
 
}{ 
 0 & \smin1 & 2 & 1 \\
 2 & 0 & 0 & 0 
 }{ 
 2 & \smin1 & 0 & \smin1 \\
 1 & 1 & 1 & 1 
 }{ 
 1 & 2 & 1 & 0 \\
 1 & \smin1 & 1 & 1 
 } \myenspace .
\end{equation*}
\begin{equation*}
(-1,2):  
\recray{
 0 & \smin1 & \smin1 & 0 \\
 1 & 1 & 0 & 0 
}{ 
 \smin1 & 1 & \smin1 & \smin1 \\
 2 & 0 & 0 & 0 
 }{ 
 \smin2 & 0 & 0 & 0 \\
 1 & 1 & 1 & 1 
}{ 
 \smin1 & \smin1 & \smin1 & \smin1 \\
 1 & \smin1 & 1 & 1 
  } \myenspace .
\end{equation*}
\begin{equation*}
\label{eqn: 2neg2}
(2,-2):  
\recray{
 1 & 1 & \smin1 & 0 \\
 \smin1 & \smin2 & \smin1 & 1 
}{ 
  2 & 0 & \smin1 & \smin1 \\
 \smin2 & 1 & 0 & \smin1 
 }{ 
0 & 2 & 1 & 1 \\
 \smin2 & \smin1 & 0 & \smin1 
 }{ 
 0 & \smin2 & 1 & 1 \\
 \smin1 & 0 & \smin1 & \smin2 
 } \myenspace .
\end{equation*}
At this point we have implied all indices of sup-norm at most $2$.  
\begin{equation*}
(3,0):  
\recray{
 2 & 1 & 0 & 0 \\
 0 & 0 & 1 & 0 
}{ 
 3 & 1 & 0 & 0 \\
 0 & 0 & 1 & 1 
 }{ 
 1 & 1 & 2 & 2 \\
 1 & \smin1 & 0 & 0 
 }{ 
 2 & \smin2 & 1 & 1 \\
 1 & 1 & 0 & 0 
 } \myenspace .
\end{equation*}
\begin{equation*}
(3,1):  
\recray{
 2 & 1 & 0 & 0 \\
 1 & 0 & 1 & 0 
}{ 
 3 & 1 & 0 & 0 \\
 1 & 1 & 1 & 1 
}{ 
 1 & 1 & 2 & 2 \\
 1 & \smin1 & 1 & 1 
}{ 
 2 & \smin2 & 1 & 1 \\
 2 & 0 & 0 & 0 
} \myenspace .
\end{equation*}
\begin{equation}\label{eqn:32}
(3,2):  
\recray{
 2 & 1 & 0 & 0 \\
 1 & 1 & 1 & 0 
}{ 
 3 & 1 & 0 & 0 \\
 2 & 0 & 1 & 1 
 }{ 
 1 & 1 & 2 & 2 \\
 2 & 0 & 1 & 1 
 }{ 
 2 & \smin2 & 1 & 1 \\
 2 & 0 & 1 & 1 
} \myenspace .
\end{equation}
\begin{equation*}
(3,3):  
\recray{
 2 & 1 & 1 & 0 \\
 2 & 1 & 0 & 0 
}{ 
 3 & 1 & 1 & 1 \\
 3 & 1 & 0 & 0 
 }{ 
 2 & 0 & 2 & 2 \\
 1 & 1 & 2 & 2 
 }{ 
 3 & \smin1 & 1 & 1 \\
 2 & \smin2 & 1 & 1 
} \myenspace .
\end{equation*}
Simply by switching top rows with bottom rows, we similarly imply $(0,3)$, $(1,3)$, and $(2,3)$.  And by putting negatives on the second row of \eqref{eqn:32}, we imply the index $(3,-2)$ (and $(-2,3)$ by switching top and bottom).
\begin{equation*}
(3,-1):  
\recray{
 2 & 1 & 0 & 0 \\
 \smin1 & \smin1 & \smin2 & 2 
}{ 
 3 & 1 & 0 & 0 \\
 \smin1 & 1 & 1 & \smin1 
 }{ 
 1 & 1 & 2 & 2 \\
 \smin1 & \smin1 & 1 & \smin1 
 }{ 
 2 & \smin2 & 1 & 1 \\
 0 & 0 & 0 & \smin2 
} \myenspace .
\end{equation*}
\begin{equation*}
(3,-3):  
\recray{
 1 & 2 & 1 & 0 \\
 \smin2 & \smin1 & 0 & 0 
}{ 
 3 & \smin1 & 1 & 1 \\
 \smin3 & \smin1 & 0 & 0 
 }{ 
 3 & 1 & 1 & 1 \\
 \smin1 & \smin1 & \smin2 & \smin2 
 }{ 
 2 & 0 & 2 & 2 \\
 \smin2 & 2 & \smin1 & \smin1 
} \myenspace .
\end{equation*}  
Again by switching top and bottom we get $(-1,3)$.  We have now implied all indices with sup-norm at most $3$.
\begin{equation*}
(4,0):  
\recray{
 2 & 1 & 0 & 1 \\
 0 & 0 & 1 & 0 
}{ 
 4 & 1 & 1 & 0 \\
 0 & 0 & 1 & 1 
 }{ 
 2 & 1 & 3 & 2 \\
 1 & \smin1 & 0 & 0 
 }{ 
 3 & \smin2 & 2 & 1 \\
 1 & 1 & 0 & 0 
  } \myenspace .
\end{equation*}\begin{equation*}
(4,1):
\recray{
 3 & 2 & 1 & \smin1 \\
 0 & 0 & 0 & 1 
}{
 4 & 1 & 0 & 1 \\
 1 & 0 & 1 & 0 
 }{ 
 2 & 1 & 2 & 3 \\
 1 & 0 & 1 & 0 
 }{ 
 3 & \smin2 & 1 & 1 \\
 1 & 0 & 1 & 0 
  } \myenspace .
\end{equation*}
\begin{equation*}
(4,2): \recray{
 3 & 2 & 1 & \smin1 \\
 1 & 1 & 1 & 0 
}{
 4 & 1 & 0 & 1 \\
 2 & 0 & 1 & 1 
}{ 
  2 & 1 & 2 & 3 \\
 2 & 0 & 1 & 1 
}{ 
  3 & \smin2 & 1 & 2 \\
 2 & 0 & 1 & 1 
} \myenspace .
\end{equation*}
\begin{equation*}
(4,3):  
\recray{
 2 & 2 & 1 & 0 \\
 2 & 1 & 0 & 0 
}{ 
 4 & 0 & 1 & 1 \\
 3 & 1 & 0 & 0 
 }{ 
 3 & 1 & 2 & 2 \\
 1 & 1 & 2 & 2 
 }{ 
 3 & \smin1 & 2 & 2 \\
 2 & \smin2 & 1 & 1 
  } \myenspace .
\end{equation*}
\begin{equation*}
(4,4):  
\recray{
 3 & 2 & 1 & \smin1 \\
 2 & 2 & 1 & 0 
}{ 
 4 & 1 & 0 & 1 \\
 4 & 0 & 1 & 1 
}{ 
 2 & 1 & 2 & 3 \\
 3 & 1 & 2 & 2 
}{ 
 3 & \smin2 & 1 & 2 \\
 3 & \smin1 & 2 & 2 
} \myenspace .
\end{equation*}
Again by switching top rows with bottom rows, we similarly imply $(0,4)$, $(1,4)$, $(2,4)$ and $(3,4)$. And by putting negatives on the second rows, we imply the indices $(4,-1)$, $(-1,4)$, $(4,-3)$ and $(-3,4)$.  
\begin{equation*}
(4,-2): \recray{
 2 & 1 & \smin1 & 1 \\
 \smin1 & \smin1 & \smin1 & 0 
}{
 4 & 1 & 0 & \smin1 \\
 \smin2 & 0 & \smin1 & \smin1 
}{ 
  1 & 2 & 3 & 2 \\
 \smin2 & 0 & \smin1 & \smin1 
}{ 
  2 & \smin3 & 2 & 1 \\
 \smin2 & 0 & \smin1 & \smin1 
} \myenspace .
\end{equation*}
\begin{equation*}
(4,-4):  
\recray{
 2 & 1 & \smin1 & 1 \\
 \smin2 & \smin2 & \smin1 & 0 
}{ 
 4 & 1 & 0 & \smin1 \\
 \smin4 & 0 & \smin1 & \smin1 
}{ 
 1 & 2 & 3 & 2 \\
 \smin3 & \smin1 & \smin2 & \smin2 
}{ 
 2 & \smin3 & 2 & 1 \\
 \smin3 & 1 & \smin2 & \smin2 
} \myenspace .
\end{equation*}
By switching rows, we imply $(-2,4)$.  We have now demonstrated the calculation of all terms of index with sup-norm at most $4$.  The second part of the statement follows immediately from the first.
\end{proof}

\subsection{Basesets for ranks $n \geq 3$}

Let $\e_i$ denote the standard basis vectors.

\begin{lemma}
\label{lem:3}
Define the followings subsets of $\Z^3$.
\begin{align*}
L_2 &= \{ \e_i \}_i \cup \{ \e_i \pm \e_j \}_{i \neq j} \cup \{ 2\e_i \}_i \\
L'_2 &= \{ a_i\e_i+a_j\e_j : a_i \in \Z, 1 \le i \le j \le 3\}.
\end{align*}
Then all indices $\v \in \Z^3$ with $N(\v)\le2$ are $L_2$-integrally implied by $L'_2$.

\end{lemma}

\begin{proof}
We make use of the recurrences
\begin{gather}
\label{eqn: recn3a}
\recray{
 1 & 1 & 0 & \smin1 \\
 0 & 0 & \smin1 & 1 \\
 1 & 0 & 1 & 0 
}{  
 1 & 0 & \smin1 & 0 \\
 1 & 0 & 0 & \smin1 \\
 1 & 1 & 1 & 1 
 }{ 
 0 & 1 & 0 & 1 \\
 0 & 1 & 1 & 0 \\
 1 & \smin1 & 1 & 1 
 }{ 
 0 & \smin1 & 0 & 1 \\
 0 & \smin1 & 1 & 0 \\
 2 & 0 & 0 & 0 
  } \myenspace , \\
\label{eqn: recn3}
\recray{
 0 & 0 & 1 &\smin1 \\
1 & 1 & 0 & \smin1 \\
 0 & 1 & 1 & 0 
}{ 
 \smin1 & 0 & 0 & 1 \\ 
 1 & 0 & \smin1 & 0 \\
 1 & \smin1 & 1 & 1 
 }{ 
 0 & \smin1 & \smin1 & 0 \\ 
 0 & 1 & 0 & 1 \\
2 & 0 & 0 & 0 
 }{ 
 0 & 1 & \smin1 & 0 \\  
 0 & \smin1 & 0 & 1 \\
 1 & 1 & 1 & 1 
  } \myenspace , \\
\label{eqn: recn3b}
\recray{
 1 & 0 & 1 & 0 \\
 0 & 0 & 0 & 1 \\
 1 & 1 & 0 & \smin1
}{ 
 1 & 1 & 1 & 1 \\
 1 & 0 & 1 & 0 \\
 1 & 0 & \smin1 & 0 
 }{
 1 & \smin1 & 1 & 1 \\
 1 & 0 & 1 & 0 \\
 0 & 1 & 0 & 1
 }{ 
 2 & 0 & 0 & 0 \\
 1 & 0 & 1 & 0 \\
 0 & \smin1 & 0 & 1 
  } \myenspace .
\end{gather}
Permute the rows of \eqref{eqn: recn3a} by the cyclic permutations $(123)$ and $(132)$, calling the results $(\ref{eqn: recn3a}')$ and $(\ref{eqn: recn3a}'')$ respectively (for example, the rightmost column of $(\ref{eqn: recn3a}')$ is $(0,1,0)$).  Do the same for \eqref{eqn: recn3} and \eqref{eqn: recn3b}.
\par
Consider the equation obtained by the combination \tiny
\begin{gather*}
(\ref{eqn: recn3a}) \times T_{(1,1,1)}T_{(1,0,0)}^2T_{(1,-1,0)}T_{(0,1,0)}^2 
+ (\ref{eqn: recn3a}') \times T_{(1,1,1)}T_{(1,0,0)}T_{(0,1,-1)}T_{(0,1,0)}^2T_{(0,0,1)} \\
+ (\ref{eqn: recn3b}) \times T_{(1,-1,0)}T_{(0,1,0)}^2T_{(0,1,1)}T_{(0,0,1)}T_{(1,0,1)}^2 
+ (\ref{eqn: recn3b}') \times T_{(0,1,-1)}T_{(1,0,0)}^2T_{(0,0,1)}T_{(1,0,1)}T_{(1,1,0)} \\
+(\ref{eqn: recn3}) \times T_{(1,1,1)}T_{(1,0,0)}^2T_{(0,1,0)}^2T_{(1,1,0)} 
+ (\ref{eqn: recn3}') \times T_{(1,1,1)}T_{(0,1,0)}^2T_{(1,0,0)}T_{(0,1,1)}T_{(0,0,1)} \\
+ (\ref{eqn: recn3}'') \times T_{(1,1,1)}T_{(1,0,0)}^2T_{(1,0,1)}T_{(0,1,0)}T_{(0,0,1)} 
\end{gather*}\normalsize
The result has the form $aT_{(1,1,1)}+b=0$ where $a$ and $b$ are polynomials in $T_\v$ where every $\v$ has at least one zero coordinate.  In particular,
\[
a = T_{(1,0,0)}^3T_{(0,1,0)}T_{(0,0,1)}^2T_{(1,0,1)}T_{(0,2,0)}T_{(1,0,-1)}.
\]
Thus $T_{(1,1,1)}$ is $L_2$-integrally implied by $L_2'$.  To imply the terms $T_{(-1,1,1)}$, $T_{(1,-1,1)}$, and $T_{(1,1,-1)}$, use $(\ref{eqn: recn3a})$, $(\ref{eqn: recn3a}')$, and $(\ref{eqn: recn3a}'')$.  This covers all terms of sup-norm at most $1$.
\par
We have the following recurrence:
\begin{equation*}
\recray{
0 & 0 & 0 & 1 \\
0 & 0 & \smin1 & 1 \\
2 & 1 & 1 & \smin1 \\
0 & 1 & 0 & 1 \\
1 & 1 & 0 & 0
}
{ 
1 & 0 & 1 & 0 \\
  1 & 0 & 0 & \smin1 \\
 2 & 1 & 0 & 1 \\
 2 & \smin1 & 1 & 0 \\
2 & 0 & 0 & 0 
  }
{ 
1 & 0 & 1 & 0 \\
  0 & 1 & 1 & 0 \\
 1 & 0 & 1 & 2 \\
 2 & 1 & 1 & 0 \\
1 & 1 & 1 & 1
  }
{ 
1 & 0 & 1 & 0 \\
  0 & \smin1 & 1 & 0 \\
 2 & \smin1 & 0 & 1 \\
 1 & 0 & 2 & 1 \\
1 & \smin1 & 1 & 1
   } \myenspace .
\end{equation*}
If $\v$ has exactly one coordinate of value $\pm2$ (the rest $\pm 1$), then we imply $\v$ by taking the first three rows in the recurrence above (possibly taking negatives and permutations of rows as necessary).  If $\v$ has exactly two $\pm 2$'s, use the middle three rows in the same way.  If $\v$ has exactly three $\pm 2$'s, use the last three rows (this relies on the previous cases).
\end{proof}
%\pagebreak

\begin{remark}
\label{remark:det}
The four equations $(\ref{eqn: recn3a})$, $(\ref{eqn: recn3a}')$, $(\ref{eqn: recn3a}'')$ and $(\ref{eqn: recn3})$ in the four unknowns $T_{(1,1,1)}$, $T_{(-1,1,1)}$, $T_{(1,-1,1)}$ and $T_{(1,1,-1)}$, are linear with coefficients consisting of monomials in $T_\v$ where $\v$ has at least one zero coordinate.  The determinant of the system is
\[
2T_{(1,0,0)}T_{(0,1,0)}T_{(0,0,1)}^2T_{(1,1,0)}T_{(1,0,1)}^2T_{(0,1,1)}^2
T_{(1,-1,0)}T_{(1,0,-1)}T_{(0,1,-1)}.
\]
This observation is useful for calculations where $2$ is invertible.
\end{remark}

\begin{theorem}
\label{thm: nlaurent}
Let $n \ge 2$.  For each $\ell$ in the set
\[
L = \{ 0, 1 \}^n \setminus \{ (0,0,\ldots, 0), (1,1,\ldots, 1) \},
\]
choose a vector $\mathbf{x}_\ell$ having $N(\mathbf{x}_\ell)=1$ and having non-zero entries exactly where $\ell$ does.  Let $G_n = \{ \mathbf{x}_\ell \}_{\ell \in L}$.  Let 
\begin{align*}
H_n &= G_n \cup \{ \e_i \} \cup \{ \e_i \pm \e_j, i \neq j\} \cup \{ 2\e_i \},\\
H'_n &= H_n \cup \{ 2\e_i + \e_j, i \neq j \}.
\end{align*}
Then $\Z^n$ is $H_n$-integrally implied by $H_n'$.
\end{theorem}

\begin{proof}
The proof is by induction on $n$.  The base case is $n=2$, which is a consequence of Theorem \ref{thm: 2laurent}.  

Fixing any $1 \le i \le n$, we can identify $H_{n-1}$ with a subset of $H_n$ (and $H_{n-1}'$ with a subset of $H_n'$) by adding a zero between the $(i-1)$-th and $i$-th positions of each vector of $H_{n-1}$ (or $H_{n-1}'$).  By this identification and by the inductive hypothesis (for $n-1$), any $\v \in \Z^n$ with a zero in the $i$-th position is $H_n$-integrally implied by $H_n'$.  Therefore it suffices to imply those $\v \in \Z^n$ having no zero coordinate.

\par
The inductive step is itself an induction on the sup-norm of $\v$.  The base cases are $N(\v)=1$ and $N(\v)=2$.  Both of these for $n=3$ are provided by Lemma \ref{lem:3}, so for the base cases, we may assume $n \ge 4$.  To imply $\v$, we construct a recurrence row-by-row, so that the first column is exactly $\v$.  For the first three rows, use the following, multiplied by $-1$ as necessary.
\begin{equation*}
\recray{
1 & 0 & 0 & 0 \\
0 & 1 & 0 & 0 \\
0 & 0 & 0 & 1 
}{ 
  1 & 1 & 0 & 0 \\
1 & \smin1 & 0 & 0 \\
 1 & 0 & 1 & 0 
  }{ 
  0 & 0 & 1 & 1 \\
1 & 1 & 0 & 0 \\
 1 & 0 & 1 & 0 
  }{ 
  1 & \smin1 & 0 & 0 \\
0 & 0 & 1 & 1 \\
 1 & 0 & 1 & 0 
} \myenspace .
\end{equation*}
For all subsequent rows, use one of the following two recurrences (shown together in an array), multiplied by $-1$ as appropriate:
\begin{equation*}
\recray{
1 & 1 & 1 & \smin1\\
0 & 0 & \smin1 & 1
}{ 
 1 & 0 & 0 & 1 \\
1 & 0 & 0 & \smin1
  }{ 
 1 & 0 & 0 & 1 \\
0 & 1 & 1 & 0
  }{ 
 1 & 0 & 0 & 1 \\
0 & \smin1 & 1 & 0 
   } \myenspace .
\end{equation*}
For each row, the choice between the two possibilities can be made in such a way that the fourth column of the recurrence lies in $G_n$.  Columns two and three have at most two non-zero entries (which are $\pm 1$) and so are in $H_n$.  The other columns (5-12) have at least one zero entry, and so are already implied by the inductive step.  This completes the case $N(\v)=1$.
\par
For the remainder of the proof, we will repeatedly use the following recurrences.  Let $w_i = \lceil \frac{v_i}{2} \rceil$.  If $v_i$ is even, we call the following recurrences (shown here in an array) $(E1)$ through $(E4)$:
\begin{equation*}
\recray{
w_i\smins1 &w_i &0 &1 \\
%w_i &w_i\splus1 & 1 & \smins1 \\
w_i & w_i\smins1 & 0 & 1 \\
%w_i\splus1 & w_i & 1 & \smins1 \\
w_i & w_i & 0 & 0 \\
w_i & w_i & 1 & 0
}{
v_i & \smins1 & 1 & 0 \\
%v_i & \smins1 & 0 & 1 \\
v_i & 1 & 1 & 0 \\
%v_i & 1 & 0 & 1 \\
v_i & 0 & 0 & 0 \\
v_i & 0 & 1 & 1
}{
w_i\splus1 & w_i & w_i & w_i\smins1 \\
%w_i \splus 1 & w_i & w_i \smins1 & w_i \\
w_i & w_i\smins1 & w_i \splus1 & w_i \\
%w_i & w_i \smins1 & w_i & w_i \splus1 \\
w_i & w_i & w_i & w_i \\
w_i \splus 1 & w_i \smins 1 & w_i & w_i 
}{
w_i & \smins w_i\splus1 & w_i \splus1 & w_i \\
%w_i & \smins w_i \splus 1 & w_i & w_i\splus1 \\
w_i\splus1 & \smins w_i & w_i & w_i\smins1 \\
%w_i\splus1 & \smins w_i & w_i\smins1 & w_i \\
w_i & \smins w_i & w_i & w_i \\
w_i\splus1 & \smins w_i\splus1 & w_i & w_i
}
\end{equation*}
If $v_i$ is odd, we call the following recurrences $(O1)$ through $(O5)$.
\begin{equation*}
\recray{
w_i & w_i\smins 1 & 0 & 0 \\
w_i \smins1 & w_i & 0 & 0 \\
%w_i & w_i\smins 1 & 1 & 0 \\
w_i \smins 1 & w_i & 1 & 0 \\
w_i & w_i & 0 & \smins1 \\
%w_i & w_i & 0 & 1 \\
w_i & w_i & 1 & \smins1 \\
%w_i & w_i & \smins 1 & 1
}{
v_i & 1 & 0 & 0 \\
v_i & \smins 1 & 0 & 0 \\
%v_i & 1 & 1 & 1 \\
v_i & \smins 1 & 1 & 1 \\
v_i & 0 & \smins 1 & 0 \\
%v_i & 0 & 1 & 0 \\
v_i & 0 & 0 & 1 \\
%v_i & 0 & 0 & \smins 1
}{
w_i \smins 1 & w_i \smins 1 & w_i & w_i \\
w_i & w_i & w_i\smins1 & w_i \smins 1 \\
%w_i & w_i\smins 2 & w_i & w_i \\
w_i\splus 1 & w_i \smins 1 & w_i \smins 1 & w_i \smins 1 \\
w_i \smins 1 & w_i & w_i \smins 1 & w_i \\
%w_i \splus 1& w_i & w_i \splus 1 & w_i \\
w_i & w_i \smins 1 & w_i \smins 1 & w_i \\
%w_i & w_i \splus 1 & w_i \splus 1 & w_i 
}{
w_i & \smins w_i & w_i \smins 1 & w_i \smins 1 \\
w_i \smins 1 & 1 \smins w_i & w_i & w_i \\
%w_i \splus 1 & \smins w_i \splus 1 & w_i \smins 1 & w_i \smins 1 \\
w_i & \smins w_i & w_i & w_i \\
w_i\smins 1 & \smins w_i & w_i \smins 1 & w_i \\
%w_i \splus 1 & \smins w_i & w_i \splus 1 & w_i \\
w_i & 1 \smins w_i &w_i \smins 1 & w_i\\
%w_i & \smins w_i \smins 1 & w_i \splus 1 & w_i
}
\end{equation*}
\par
The second base case is $N(\v)=2$ ($n \ge 4$ still).  Since we may assume $v_i \neq 0$ (this is covered by previous cases in the induction on $n$), the other $v_i$ have $|v_i| = \pm 1$.  There are three cases:

{\bf Case I: $\v$ has at least three odd $v_i$}.  Use for the first three odd $v_i$ the recurrences $(O1)$, $(O4)$ and $(O5)$ respectively.  Use $(E3)$ for all the even $v_i$.  In this case, all the columns besides the first contain only digits $0$ and $\pm1$ and so were implied in the case $N(\v) = 1$.  The columns 2-4 contain only one non-zero term each, and so are in $H_n$.  

{\bf Case II: $\v$ has one or two odd $v_i$}.  Use $(O3)$ for one odd coordinate and $(O1)$ for the other (if it exists).  Use $(E3)$ for all even coordinates.  Then, the columns 2-4 contain one or two non-zero entries, and the columns 5-12 may contain at most one $\pm 2$, but such a column was implied in the Case I.

{\bf Case III: $\v$ has no odd $v_i$}.  Use $(E1)$ and $(E4)$ for the first two rows, and $(E3)$ for all others.  Columns 2-4 contain one or two non-zero entries and 5-12 at most two $\pm 2$'s, but such a column was implied in Case I or II.

This completes the $N(\v)=2$ base case.
\par
Now suppose  $N(\v) = N_0 \ge 3$ and $n \ge 3$.  This is the inductive step; we will assume we have implied all indices of sup-norm less than $N_0$.  As before, $v_i \neq 0$.  
For $|v_i|=3$, $(O1)$, $(O2)$, $(O4)$, and $(O5)$ have entries less than $N_0$ in columns 5-12.  For $1 \le |v_i| \le 2$, and $3 < |v_i| \le N_0$, all applicable recurrences have entries less than $N_0$ in those columns.  We have two cases:

{\bf Case I: $\v$ has at least one even entry}.  Use $(E4)$ for the first even coordinate, and choose from $(E1)$ and $(E2)$ for the second even coordinate (if it exists).  We use $(E3)$ for all other even coordinates.  We will use $(O1)$ or $(O2)$ for all odd entries (and make the choice between $(E1)$ and $(E2)$ above) in such a way that the second column is in $G_n$.  

{\bf Case II:  $\v$ has no even entry}.  Use $(O4)$ and $(O5)$ for the first two odd coordinates, and $(O1)$ or $(O2)$ for all others, according so that the second column is an element of $G_n$.
\end{proof}

\section{Net polynomials over $\C$}
\label{sec: overc}

Fix an elliptic curve $E$ defined over $\C$.  Our purpose is to define rational functions $\Omega_\v: E^n \rightarrow \C$ for all $\v \in \Z^n$ such that for each $\mathbf{P} \in E^n$, the map
\[
W_{E,\mathbf{P}}: \Z^n \rightarrow \C, \quad \v \mapsto \Omega_\v(\mathbf{P})
\]
is an elliptic net.  In this section we associate a lattice $\Lambda \subset \C$ to the elliptic curve $E$ and consider the complex uniformization $\C / \Lambda$.

\subsection{Elliptic functions over $\C$}

For a complex lattice $\Lambda$, let $\eta : \Lambda \rightarrow \C$ be the quasi-period homomorphism, and define a quadratic form $\lambda : \Lambda \rightarrow \{ \pm 1 \}$ by 
\begin{equation*}
\lambda(\omega) = \left\{ 
\begin{array}{ll}
1 & \text{if }\omega \in 2\Lambda\text{,}  \\
-1 & \text{if }\omega \notin 2\Lambda\text{.}
\end{array} \right.
\end{equation*}
Recall that the Weierstrass sigma function $\sigma: \C / \Lambda \rightarrow \C$ satisfies the following transformation formula for all $z \in \C$ and $\omega \in \Lambda$:
\begin{equation}
\label{eqn: sigma}
\sigma(z+\omega; \Lambda) = \lambda(\omega) e^{\eta(\omega)(z+\frac{1}{2}\omega)}\sigma(z;\Lambda)
\end{equation}
\begin{definition}
\label{defn: psi}
Fix a lattice $\Lambda \in \C$ corresponding to an elliptic curve $E$.  For $\mathbf{v} = (v_1,\ldots,v_n) \in \Z^n$, define a function $\Omega_{\mathbf{v}}$ on $\C^n$ in variables $\mathbf{z} = (z_1,\ldots,z_n)$ as follows:
\[ \Omega_{\mathbf{v}}(\mathbf{z};\Lambda) = \frac{\sigma(v_1z_1 + \ldots + v_nz_n;\Lambda)}{\displaystyle{\prod_{i=1}^n \sigma(z_i;\Lambda)^{2v_i^2-\sum_{j=1}^nv_iv_j}\prod_{1 \leq i < j \leq n }\sigma(z_i+z_j;\Lambda)^{v_iv_j}}}. \]
(If $\v = \mathbf{0}$, we set $\Omega_\v \equiv 0$.)
In particular, we have for each $n\in\Z$, a function $\Omega_{n}$ on $\C$ in the variable $z$:
\[ \Omega_{n}(z;\Lambda) = \frac{\sigma(nz;\Lambda)}{\sigma(z;\Lambda)^{n^2}}, \]
and for each pair $(m,n) \in \Z\times \Z$, a function $\Omega_{m,n}$ on $\C \times \C$ in variables $z$ and $w$:
\[ \Omega_{m,n}(z,w;\Lambda) = \frac{\sigma(mz + nw;\Lambda)}{\sigma(z;\Lambda)^{m^2-mn}\sigma(z+w;\Lambda)^{mn}\sigma(w;\Lambda)^{n^2-mn}}. \]
\end{definition}

\begin{remark} Compare the proof of Lemma \ref{lem:quad} to this definition. \end{remark}

\begin{proposition}
\label{prop: psielliptic}
Fix a lattice $\Lambda \in \C$ corresponding to an elliptic curve $E$.
The functions $\Omega_{\mathbf{v}}$ are elliptic functions in each variable.
\end{proposition}

\begin{proof}
Let $\omega \in \Lambda$.  We show the function is elliptic in the first variable.  Let $\mathbf{v} = (v_1,\ldots, v_n) \in \Z^n$ and $\mathbf{z} = (z_1,\ldots,z_n), \mathbf{w} = (\omega, 0, \ldots, 0) \in \C^n$.  Using \eqref{eqn: sigma}, we calculate
$$\frac{\Omega_\mathbf{v}(\mathbf{z}+\mathbf{w}; \Lambda)}{\Omega_\mathbf{v}(\mathbf{z}; \Lambda)} 
= \frac{\lambda(v_1\omega)}{\lambda(\omega)^{v_1^2}} = 1$$
where the last equality holds because $\lambda$ is a quadratic form.
%If $\omega, v_1\omega \notin 2\Lambda$, then $v_1$ is odd, and $F = 1$.  If $\omega \notin 2\Lambda$ but $v_1\omega \in 2\Lambda$, then $v_1$ must be even, and so $F = 1$ again.  Finally, if $\omega \in 2\Lambda$, then $v_1\omega \in 2\Lambda$, and $F=1$.
Thus $\Omega_{\mathbf{v}}$ is invariant under adding a period to the variable $z_1$.  Similarly $\Omega_{\mathbf{v}}$ is elliptic in each variable on $(\C/\Lambda)^n$. 
\end{proof}

\begin{proposition}
\label{prop: psiuu}
Fix a lattice $\Lambda \in \C$.  Let $\mathbf{v} \in \Z^m$ and $\mathbf{z} \in \C^n$.  Let $T$ be an $n \times m$ matrix with entries in $\Z$ and transpose $T^{tr}$.  Then
\[ \Omega_{\mathbf{v}}(T^{tr}(\mathbf{z});\Lambda) = \frac{
\Omega_{T(\mathbf{v})}(\mathbf{z};\Lambda)
}{
\displaystyle{
\prod_{i=1}^n \Omega_{T(\mathbf{e}_i)}(\mathbf{z};\Lambda)^{2v_i^2-\sum_{j=1}^nv_iv_j}
\prod_{1 \leq i < j\leq n}\Omega_{T(\mathbf{e}_i +\mathbf{e}_j)}(\mathbf{z};\Lambda)^{v_iv_j}
}
} \text{.} \]
\end{proposition}

\begin{proof} A straightforward calculation using Definition \ref{defn: psi}.
\end{proof}

Let $\wp$ and $\zeta$ denote the usual Weierstrass functions.

\begin{lemma}
\label{lemma: ellids}
%Whenever $u, v \notin \Lambda$,
\[
\wp(u) - \wp(v) = - \frac{\sigma(u+v)\sigma(u-v)}{\sigma(u)^2\sigma(v)^2} \myenspace .
\]
%Furthermore, whenever $\v \cdot \z, \w \cdot \z \notin \Lambda$,
\[
\wp(\v \cdot \z) - \wp( \w \cdot \z) = - \frac{\Omega_{\v + \w}(\z)\Omega_{\v - \w}(\z)}{\Omega_\v(\z)^2\Omega_\w(\z)^2} \myenspace .
\]
\end{lemma}

\begin{proof}
The first statement is well-known (e.g. \cite{Cha}).  The second statement follows by direct calculation using Definition \ref{defn: psi}.
\end{proof}

\begin{lemma}
\label{lemma: zetaid}
%Whenever $a, b, x+a, x+b \notin \Lambda$,
\begin{gather*}
\label{eqn: zetarec}
\zeta(x+a) - \zeta(a) - \zeta(x + b) + \zeta(b) \notag \\= \frac{\sigma(x+a+b)\sigma(x)\sigma(a-b)}{\sigma(x+a)\sigma(x+b)\sigma(a)\sigma(b)} \myenspace ,
\end{gather*}
%Furthermore, whenever $x, x+a, x+b, x+a+b \notin \Lambda$,
\begin{gather*}
\zeta(x+a+b) - \zeta(x+a) - \zeta(x+b) + \zeta(x) \notag \\=  \frac{\sigma(2x+a+b)\sigma(a)\sigma(b)}{\sigma(x+a+b)\sigma(x+a)\sigma(x+b)\sigma(x)} \notag \myenspace .
\end{gather*}
\end{lemma}

\begin{proof}
Denote by $f$ and $g$ the left and right side of the first equation respectively.  Considered as functions of any one of $x$, $a$ or $b$, these are elliptic functions.  Suppose that $a, b \notin \Lambda$.  Consider $f$ and $g$ as functions of $x$.  The set of poles of $f$ or $g$ is $\{ -a, -b \}$.  The zeroes of $g$ are at $-a-b$ and $0$.  These are also zeroes of $f$, since $\zeta$ is an odd function.  Hence $f = cg$ for some $c$ not depending on $x$.  Now define instead
\begin{align*}
F &= \left( \zeta(x+a) - \zeta(a) - \zeta(x + b) + \zeta(b) \right) \sigma(x+a)\sigma(x+b) \myenspace , \\
G &= \sigma(x+a+b)\sigma(x) \myenspace .
\end{align*}
We have $F = c'G$ for some constant $c'$ independent of $x$.  Taking derivatives and evaluating at $x=0$, we have
\begin{equation*}
\left( \wp(b) -\wp(a) \right) \sigma(a)\sigma(b) = c' \sigma(a+b)\sigma'(0) 
\end{equation*}
We have $\sigma'(0) = 1$.  By Lemma \ref{lemma: ellids}, we then have
$$ c' = - \frac{\sigma(a-b)}{\sigma(a)\sigma(b)} $$
which proves the first equation.  The second is obtained by a change of variables $x \leftarrow a$, $a \leftarrow x+b$, $b \leftarrow x$.
\end{proof}

\subsection{Forming the elliptic net}

\begin{theorem}
\label{thm: edsfromcurvesc}
Fix a lattice $\Lambda \in \C$ corresponding to an elliptic curve $E$.  Fix $z_1, \ldots, z_n \in \C$.  Then the function $W: \Z^n \rightarrow \C$ defined by 
$$W(\v) = \Omega_{\v}(z_1,\ldots,z_n;\Lambda)$$
is an elliptic net.
\end{theorem}

\begin{proof}
For notational simplicity, we drop the arguments $z_i, \Lambda$ on $\Omega_\v$ and also write $\sigma(\v)$, $\wp(\v)$ and $\zeta(\v)$ for $\sigma(v_1z_1 + \ldots + v_nz_n)$, $\wp(v_1z_1 + \ldots + v_nz_n)$ and $\zeta(v_1z_1 + \ldots + v_nz_n)$.  We observe that $\v = \mathbf{0}$ if and only if $\Omega_{\v} \equiv 0$.  

We intend to show that \eqref{eqn: ellrec} holds for $W$ in $\p$, $\q$, $\r$ and $\s$.  If any one of $\p$, $\q$ or $\r$ are zero, then \eqref{eqn: ellrec} holds trivially (note that $\sigma$ is an odd function, so that $\Omega_{-\v}= - \Omega_\v$).  Hence we may assume that none of $\Omega_\p$, $\Omega_\q$, or $\Omega_\r$ is identically zero.  For any quadratic form $f$ defined on $\Z^n$, we have the following relation for all $\p, \q, \s \in \Z^n$: 
\begin{equation}
\label{eqn: quadformpsi}
f( \p + \q + \s) + f(\p - \q) + f(\s) - f(\p+\s) - f(\p) - f(\q+\s) - f(\q) = 0 \myenspace .
\end{equation}
First we address the case that  $\s = \mathbf{0}$.  By \eqref{eqn: quadformpsi} and Lemma \ref{lemma: ellids},
\[
 \frac{ \Omega_{\p + \q } \Omega_{\p - \q} }{ \Omega_{\p}^2 \Omega_{\q}^2  }
= \frac{ \sigma(\p + \q) \sigma(\p - \q) }{ \sigma(\p )^2 \sigma(\q )^2 } 
= \wp(\q) - \wp(\p) \myenspace .
\]
Therefore, we have
$$
\frac{ \Omega_{\p + \q } \Omega_{\p - \q}  }{ \Omega_{\p }^2  \Omega_{\q }^2  }
+ \frac{ \Omega_{\q + \r } \Omega_{\q - \r} }{ \Omega_{\q }^2  \Omega_{\r }^2  }
+ \frac{ \Omega_{\r + \p } \Omega_{\r - \p}  }{ \Omega_{\r }^2  \Omega_{\p }^2  }
= 0 \myenspace ,
$$
which gives the relation \eqref{eqn: ellrec} for $\s = \mathbf{0}$, that is,
$$
 \Omega_{\p + \q } \Omega_{\p - \q}   \Omega_{\r }^2 
+ \Omega_{\q + \r } \Omega_{\q - \r}  \Omega_{\p }^2 
+  \Omega_{\r + \p } \Omega_{\r - \p}  \Omega_{\q }^2  
= 0 \myenspace .
$$

Now suppose that $\s \neq \mathbf{0}$ and so $\Omega_\s \not\equiv 0$.  By \eqref{eqn: quadformpsi} and Lemma \ref{lemma: zetaid},
\begin{align*}
 \frac{ \Omega_{\p + \q + \s} \Omega_{\p - \q} \Omega_{\s} }{ \Omega_{\p + \s} \Omega_{\p} \Omega_{\q + \s} \Omega_{\q} }
&= \frac{ \sigma(\p + \q + \s) \sigma(\p - \q) \sigma(\s) }{ \sigma(\p + \s) \sigma(\p) \sigma(\q + \s) \sigma(\q) } \\
&= \zeta(\p + \s) - \zeta(\p) - \zeta(\q + \s) + \zeta(\q) \myenspace .
\end{align*}
Therefore, we have
$$
\frac{ \Omega_{\p + \q + \s} \Omega_{\p - \q} \Omega_{\s} }{ \Omega_{\p + \s} \Omega_{\p} \Omega_{\q + \s} \Omega_{\q} }
+ \frac{ \Omega_{\q + \r + \s} \Omega_{\q - \r} \Omega_{\s} }{ \Omega_{\q + \s} \Omega_{\q} \Omega_{\r + \s} \Omega_{\r} }
+ \frac{ \Omega_{\r + \p + \s} \Omega_{\r - \p} \Omega_{\s} }{ \Omega_{\r + \s} \Omega_{\r} \Omega_{\p + \s} \Omega_{\p} }
= 0 \myenspace ,
$$
or, more simply,
$$
 \Omega_{\p + \q + \s} \Omega_{\p - \q}  \Omega_{\r + \s} \Omega_{\r} 
+  \Omega_{\q + \r + \s} \Omega_{\q - \r} \Omega_{\p + \s} \Omega_{\p} 
+  \Omega_{\r + \p + \s} \Omega_{\r - \p} \Omega_{\q + \s} \Omega_{\q} 
= 0 \myenspace ,
$$
which is what was required to prove.
\end{proof}

The identity \eqref{eqn: ellrec} for $\Omega_\v$ is similar to several identities known in complex function theory \cite{GasRah, WenEkhCha}.

\subsection{Explicit rational functions}

Elliptic functions for a lattice $\Lambda$ of $\C$ give rational functions on the associated elliptic curve (via complex uniformization).  If we give a Weierstrass model for the same elliptic curve, we can give explicit expressions for the rational functions as elements of the usual field of rational functions associated to the model.  In the following proposition, we do this for $\Omega_\v$ for some small $\v \in \Z^n$, for $n=1,2,3$.

\begin{proposition}
\label{prop: 2explicit}
Consider an elliptic curve $E$, and a Weierstrass model for $E$ given by
\begin{equation*}
y^2 + a_1xy + a_3y - x^3 - a_2x^2 - a_4x - a_6 =0\myenspace .
\end{equation*}
As usual, let
\begin{gather*}
b_2 = a_1^2 + 4a_2 , \quad b_4 = 2 a_4 + a_1a_3, \quad b_6 = a_3^2 + 4a_6, \\
b_8 = a_1^2a_6 + 4a_2a_6 - a_1a_3a_4 + a_2a_3^2 - a_4^2 .
\end{gather*}
To $E$ we can also associate a complex uniformization and elliptic functions $\Omega_{\mathbf v}$ as above.  As rational functions on $E$, we have the following equalities.
\par 
For $n=1$:
\begin{align*}
\Omega_1 &= 1 \myenspace , \qquad \Omega_2 = 2y + a_1x + a_3 \myenspace , \\
\Omega_3 &= 3x^4 + b_2x^3 + 3b_4x^2 + 3b_6x + b_8 \myenspace,
\end{align*}
\vspace{-2.1em}
\begin{multline*}
\Omega_4 = (2y+a_1x+a_3) \times \\
 (2x^6 +b_2x^5+5b_4x^4+10b_6x^3+10b_8x^2+(b_2b_8-b_4b_6)x + b_4b_8-b_6^2) \myenspace ;
\end{multline*}
\par
For $n=2$:
\[
\Omega_{(1,0)}=\Omega_{(0,1)}=\Omega_{(1,1)}=1,
\]
\[
\Omega_{(1,-1)} = x_2 - x_1, \qquad \Omega_{(-1,1)} = x_1-x_2,
\]
\begin{align*}
\Omega_{(2,1)} &= 2x_1 + x_2 - \left(\frac{y_2-y_1}{x_2-x_1}\right)^2 - a_1 \left(\frac{y_2-y_1}{x_2-x_1}\right) + a_2 \myenspace ,\\
\Omega_{(1,2)} &= x_1 + 2x_2 - \left(\frac{y_2-y_1}{x_2-x_1}\right)^2 - a_1 \left(\frac{y_2-y_1}{x_2-x_1}\right) + a_2 \myenspace ;
\end{align*}
\par
For $n=3$:
%\label{prop: explicit}
\[ \Omega_{(1,0,0)}=\Omega_{(0,1,0)}=\Omega_{(0,0,1)}=\Omega_{(1,1,0)}=\Omega_{(0,1,1)}=\Omega_{(1,0,1)}=1, \]
\begin{align*}
\Omega_{(1,-1,0)}&= x_2 - x_1, \;\;\;\; \Omega_{(0,1,-1)}= x_3 - x_2, \;\;\;\; \Omega_{(-1,0,1)}= x_1 - x_3,\\
\Omega_{(-1,1,0)}&= x_1 - x_2, \;\;\;\; \Omega_{(0,-1,1)}= x_2 - x_3, \;\;\;\; \Omega_{(1,0,-1)}= x_3 - x_1,
\end{align*}
\begin{align*}
\Omega_{(1,1,1)} &= \frac{ y_1 (x_2 - x_3) + y_2 (x_3 - x_1) + y_3 (x_1 - x_2) }{ (x_1 - x_2) (x_1 - x_3) (x_2 - x_3) }, \\
\Omega_{(-1,1,1)} &= \frac{ y_1 (x_2 - x_3) - y_2 (x_3 - x_1) - y_3 (x_1 - x_2) }{ (x_2 - x_3) } + a_1 x_1 + a_3\myenspace ,\\
\Omega_{(1,-1,1)} &= \frac{ - y_1 (x_2 - x_3) + y_2 (x_3 - x_1) - y_3 (x_1 - x_2) }{ (x_3 - x_1)  } + a_1 x_2 + a_3\myenspace ,\\
\Omega_{(1,1,-1)} &= \frac{ - y_1 (x_2 - x_3) - y_2 (x_3 - x_1) + y_3 (x_1 - x_2) }{ (x_1 - x_2)  } + a_1 x_3 + a_3\myenspace .
\end{align*}
\end{proposition}

\begin{proof}The division polynomial formul{\ae} (the $n=1$ case) are well-known 
\cite{Cha} \cite[p.80]{FreLan} \cite[Exercise 3.7]{Sil1}.  The formul{\ae}  for $n=2$ and the related first dozen formul{\ae} for $n=3$ are immediate consequences of Lemma \ref{lemma: ellids} and the addition law for elliptic curves \cite[Algorithm 2.3]{Sil1}.  Only the cases where $n=3$, $v_i \neq 0$ for all $i=1,2,3$ are not immediate:  these formul{\ae} are a result of the proof of Lemma \ref{lem:3}.  Note that using Remark \ref{remark:det} results in the same formul{\ae}.
\end{proof}

\section{Net polynomials over arbitrary fields}
\label{sec: netpoly}

In the last section, we defined elliptic functions $\Omega_\v$ in the case of $\C/\Lambda$.  In this section we wish to define the same rational functions for any elliptic curve over any field, calling them $\Psi_\v$, the \emph{net polynomials}.  We will start from the results of the last section.

\subsection{Defining net polynomials}

Let $R = \Q[\alpha_1,\alpha_2, \alpha_3,\alpha_4,\alpha_6]$ be a polynomial ring over $\Q$ in the variables $\alpha_i$.  Define $f(x,y) \in R[x,y]$ by
\[
f(x,y) = y^2 + \alpha_1xy + \alpha_3y - x^3 - \alpha_2x^2 - \alpha_4x - \alpha_6.
\]
Consider the affine scheme $\mathcal{E} : f(x,y) = 0$ over $R$.  Let $\mathbf a = (a_i) \in \C^5$.  The association $(\alpha_i) \mapsto (a_i)$ gives a map $\phi_\mathbf{a}: R \rightarrow \C$.  Consider the affine variety over $\C$ given by
\[
C_{\mathbf a} : y^2 + a_1xy + a_3 y = x^3 + a_2 x^2 + a_4 x + a_6.
\]
Then $\phi_{\mathbf a}$ gives rise to a Cartesian diagram
\[
\xymatrix{
\mathcal{E}^n \ar[d] & C_{\mathbf a}^n \ar[l] \ar[d] \\
\Spec(R) & \Spec(\C) \ar[l]
}
\]
where $\mathcal{E}^n=\mathcal{E}\times_{\Spec R}\cdots\times_{\Spec R}\mathcal{E}$ is the $n$-fold fibre product of $\mathcal{E}$ with itself over $R$.

The rational functions $\Omega_\v \in \mathcal{K}(C_{\mathbf{a}}^n)$ have rational expressions in $x$, $y$ and the $a_i$ (in terms of the Weierstrass model, as in for example Proposition \ref{prop: 2explicit}).  These expressions have rational coefficients by construction and the general theory of sigma functions (the divisors are Galois invariant).  So these same expressions (with $a_i$ replaced with $\alpha_i$) give rational functions $\Psi_\v \in \mathcal{K}(\mathcal{E}^n)$.

We obtain the following theorem.

\begin{theorem}
\label{thm: omegav}
Let $n \ge 1$.  Denote by $\mathcal{K}(\mathcal{E}^n)$ the field of rational functions on $\mathcal{E}^n$.  There exists a unique system of functions $\Psi_{\mathbf v} \in \mathcal{K}(\mathcal{E}^n)$ depending on $\v \in \Z^n$ such that
\begin{enumerate}
\item the map
\[
W: \Z^n \rightarrow \mathcal{K}(\mathcal{E}^n), \quad \v \mapsto \Psi_\v
\]
is an elliptic net, and
\item 
whenever $C_\mathbf{a}$ is elliptic, the restriction of $\Psi_\v$ to a fibre $C_{\mathbf a}^n$ is the rational function $\Omega_\v$.% on $C_{\mathbf a}^n$.
\end{enumerate}
\end{theorem}

\begin{proof}
The union of the $C_\mathbf{a}^n$ for which $C_\mathbf{a}$ is an elliptic curve is Zariski dense, and so the $\Psi_\v$ are determined uniquely by their restrictions to these fibres.
\end{proof}

We call these $\Psi_\v$ the \emph{net polynomials}.  Look ahead to Section \ref{sec:polyprime} for the definition of the `polynomial' ring $\mathcal{R}_n$ in which they live.

We transfer some useful properties of the $\Omega_\v$ to properties of the $\Psi_\v$ on $\mathcal{E}^n$.  Again, there are unique rational functions $X$ and $Y$ for $\mathcal{E}$ whose restriction to elliptic $C_\mathbf{a}$ correspond to the Weierstrass functions $\wp$ and $\frac{1}{2}\wp'$.  Each $\v \in \Z^n$ gives rise to a map $\v: \mathcal{E}^n \rightarrow \mathcal{E}$ which is the linear combination associated to the vector $\v$ (e.g., $(1,1)$ is the usual group law).  Define rational functions $X_\v = X \circ \v $ and $Y_\v = Y \circ \v$ on $\mathcal{E}^n$.

The next lemma follows immediately from Lemma \ref{lemma: ellids}.

\begin{lemma}
\label{lemma: ellidsgen}
\begin{equation}
\label{eqn: ellidsgen}
\Psi_\v^2\Psi_\w^2(X_\v - X_\w) = - \Psi_{\v+\w}\Psi_{\v-\w} \myenspace .\\
\end{equation}
\end{lemma}

More generally, there is a map $T: \mathcal{E}^m \rightarrow \mathcal{E}^n$ associated to any $T \in M_{n \times m}(\Z)$.  The next proposition follows from Proposition \ref{prop: psiuu}.

\begin{proposition}
\label{prop: psiuuomega}
Let $\mathbf{v} \in \Z^n$.  Let $T$ be any $n \times m$ matrix with entries in $\Z$ and transpose $T^{tr}$.  Then
\[ 
(\Psi_{\mathbf{v}} \circ T) \displaystyle{
\prod_{i=1}^n \Psi_{T^{tr}(\mathbf{e}_i)}^{2v_i^2-\sum_{j=1}^nv_iv_j}
\prod_{1 \leq i<j\leq n }\Psi_{T^{tr}(\mathbf{e}_i +\mathbf{e}_j)}^{v_iv_j} 
=
\Psi_{T^{tr}(\mathbf{v})}
} \text{.} 
\]
\end{proposition}

\subsection{Net polynomials at primes}
\label{sec:polyprime}

In this section we determine a little more about the exact nature of the elliptic net $\Psi_\v$.  In particular, we wish to restrict the possible divisor of $\Psi_\v$, and show that it has zero valuation for certain primes.

Consider the ring $S = \Z[\alpha_1,\alpha_2, \alpha_3, \alpha_4, \alpha_6]$.  Since $f(x,y)$ is defined over $S$, we may define $\mathcal{E}_S : f(x,y) = 0$ as a scheme over $\Spec S$ whose fibre over $\Spec R$ is $\mathcal{E}$.  Then $\mathcal{E}_S^n = \mathcal{E}_S \times_{\Spec S} \cdots \times_{\Spec S} \mathcal{E}_S$ is a scheme over $\Spec S$ whose fibre over $\Spec R$ is $\mathcal{E}^n$.  Define
\[
\mathcal{R}_n =  \left. S[x_i,y_i]_{1 \leq i \leq n}\left[(x_i - x_j)^{-1}\right]_{1 \leq i < j \leq n} \right/ \left< f(x_i,y_i) \right>_{1\leq i \leq n} \myenspace .
\]
The ring $\mathcal{R}_n$ is the affine coordinate ring of the affine piece of $\mathcal{E}_S^n$ obtained by removing all the diagonals and antidiagonals, in the sense of the elliptic curve group law (in other words, on an elliptic curve fibre, $x_i = x_j$ if and only if the corresponding points satisfy $P_i = \pm P_j$).  There is a natural identification of $\mathcal{R}_n$ with a subset of $\mathcal{K}(\mathcal{E}^n)$.

\begin{theorem}
\label{thm: valuation}
The functions $\Psi_\v$ are elements of $\mathcal{R}_n$.  Let $\mathfrak{p}$ be any prime of $\mathcal{R}_n$ which is a lift of a prime of $S$.  Then $\Psi_\v \notin \mathfrak{p}$.
\end{theorem}

The lifted ideal $\mathfrak{p}=\mathfrak{q}\mathcal{R}_n$ is prime whenever $\mathfrak{q}$ is a prime of $S$.  The proof of the theorem will involve showing for all valuations $v$ associated to such primes $\mathfrak{p}$ that $v(\Psi_\v)$ (slightly modified) is a quadratic form with certain vanishing.  Then the following lemma will establish that this function is identically zero.

Let $B$ and $C$ be abelian groups written additively.  The function $f: B \rightarrow C$ is a \emph{quadratic form} if for all $x,y,z \in B$,
$$f(x+y+z) - f(x+y) - f(y+z) - f(x+z) + f(x) + f(y) + f(z) = 0 \myenspace .$$
If $f$ is a quadratic form, then for all $x, y \in B$, 
\[
f( x + y ) + f(x-y) - 2f(x) - 2f(y) = 0.
\]
The converse holds if $C$ is is $2$-torsion free.

\begin{lemma}
\label{lem:quad}
Let $M: \Z^n \rightarrow \Z$ be a quadratic form.  Suppose that $M(\v)=0$ for all $\v = \e_i$ and $\v = \e_i + \e_j$ (i.e. for standard basis vectors and their two-term sums).  Then $M(\v)=0$ for all $\v$.
\end{lemma}

\begin{proof}
It is well-known that any value of a quadratic form can be given in terms of its value at a certain `base' of vectors.  In particular,
\[
f\left( \sum_{i=1}^n a_i \e_i \right) = \sum_{i=1}^n \left( 2a_i^2 - \sum_{j=1}^n a_ia_j \right)  f(\e_i) + \sum_{1 \le i < j \le n } a_ia_j f(\e_i + \e_j).
\]
\end{proof}

\begin{proof}[Proof of Theorem \ref{thm: valuation}]
Each $\Psi_\v \in \mathcal{K}(\mathcal{E}^n)$ has a corresponding Weil divisor.  Suppose a codimension-one subscheme $X$ appears as a summand in this divisor, and let $\widetilde X = X \cap C_{\mathbf a}^n$.  If $C_{\mathbf a}$ is elliptic, $\widetilde X \neq \emptyset$, and $\widetilde X \neq C_\mathbf{a}^n$, then $\widetilde X$ is of codimension one in $C_{\mathbf a}^n$ and appears in the divisor of $\Omega_\v$ to the same order as $X$ appears in the divisor of $\Psi_\v$.  Definition \ref{defn: psi} determines the divisors of $\Omega_\v$ and this restricts the possible divisors for $\Psi_\v$.  In particular, it shows that $s \Psi_\v \in \mathcal{R}_n$, where $s \in S$. 

%In particular, the codimension one subschemes of $\mathcal{E}^n$ which intersect trivially with fibres $C_{\mathbf a}^n$ are fibres of the map $\mathcal{E}^n \rightarrow \Spec R$.  
%This tells us that $\Psi_\v = r/s$ for $s \in S$ and $r \in \mathcal{R}_n$.  For the first statement of the theorem, it suffices to show that $s$ can be taken to be $1$.  To this end, let $v$ denote any valuation of $S$.  

Therefore, taking $v$ to be a valuation of $\mathcal{R}_n$ lifted from a valuation of $S$ associated to a prime $\mathfrak{q}$ of $S$, it will suffice to show that $v(\Psi_\v)=0$ for all $\v \in \Z^n$.
 
Equation \eqref{eqn: ellidsgen} of Lemma \ref{lemma: ellidsgen} implies
\[
X_\v - X_\w = - \frac{ \Psi_{\v+\w}\Psi_{\v-\w} }
 { \Psi_\v^2\Psi_\w^2 }.
\]
We claim that $v(X_\v - X_\w)=0$ whenever $\v \neq \pm \w$, $\v \neq 0$, and $\w \neq 0$.  
\par 
First we show that $v(X_\v - X_\w) < 0 \implies \v = 0 \mbox{ or } \w = 0$.
Suppose $v(X_\v - X_\w) < 0$.  Then, $v(X_\v)<0$ or $v(X_\w)<0$.  Suppose $v(X_\v) < 0$.  This implies that $\v(\mathbf{P}) = \mathcal{O}$ for all $\mathbf{P}$ on the non-singular part of the fibre over $\mathfrak{q}$ of $\mathcal{E}_S$.
%$$y^2 + \alpha_1xy + \alpha_3y = x^3 + \alpha_2x^2 + \alpha_4x + \alpha_6$$ over $S$.
Since $\mathbf{P}$ ranges over all possible values (e.g. $\mathbf{P} = (P, \mathcal{O}, \ldots, \mathcal{O})$), we find that this implies that $[v_i]=[0]$ for all $i$.  In turn, this shows that $\v = 0$.  Similarly, if $v(X_\w)<0$, then $\w = 0$.
\par
Now we show that $v(X_\v - X_\w) > 0 \implies \v = \pm \w$.
Suppose the valuation is positive.  Then $\v(\mathbf{P}) = \pm \w(\mathbf{P})$ for all $\mathbf{P}$ on the non-singular part of the fibre over $\mathfrak{q}$ of $\mathcal{E}_S$.  Since $\mathbf{P}$ ranges over all possible values (e.g. $\mathbf{P} = (P, \mathcal{O}, \ldots, \mathcal{O})$ or $\mathbf{P}=(P, P, \mathcal{O}, \ldots, \mathcal{O})$), we find that this implies, in particular, that for all $0 \le i \le j \le n$, we have $[v_i \pm w_i]=[0]$ and $[v_i+v_j \pm (w_i + w_j)] = [0]$ on $\mathcal{E}_S$.  In turn, this gives $v_i = \pm w_i$ and $v_i + v_j = \pm (w_i+w_j)$.  Together these imply that $\v = \pm \w$.  This demonstrates the claim.
\par
Define a function $M: \Z^n \rightarrow \Z$ by
\[
M(\v) = \left\{ \begin{array}{ll}
v(\Psi_\v) & \v \neq 0 \\
0 & \v = 0
\end{array}
\right.
\]
Note that $M(-\v) = M(\v)$, from which one can deduce that
\begin{equation}
\label{eqn:parallel}
M(\v+\w) + M(\v-\w) - 2M(\v) - 2M(\w) = 0
\end{equation}
whenever $\v=0$ or $\w=0$.  Our work up until now has shown that \eqref{eqn:parallel} holds in all other cases except $\v+\w=0$ or $\v-\w=0$.  These remaining two cases reduce to the statement that for all $\mathbf{u}$, $M(2\mathbf{u}) = 4M(\mathbf{u})$.  To obtain this, take the sum of the four instances of \eqref{eqn:parallel} with $(\v,\w)$ respectively taking the values $(4\mathbf{u}, \mathbf{u})$, $(3\mathbf{u}, \mathbf{u})$, $(3\mathbf{u}, \mathbf{u})$ and $(2\mathbf{u}, \mathbf{u})$, and then subtract the instance of \eqref{eqn:parallel} with $(\v,\w) = (3\mathbf{u}, 2\mathbf{u})$.  
\par
We have shown that \eqref{eqn:parallel} holds for all $\v$ and $\w$, and that therefore $M: \Z^n \rightarrow \Z$ is a quadratic form (since $\Z$ is $2$-torsion free).  The other assumptions of Lemma \ref{lem:quad} are verified by Proposition \ref{prop: 2explicit}.  Therefore, $M$ is identically zero, which is what was required to prove.
\end{proof}

\subsection{Summary}

Let $n \geq 1$.  For any elliptic curve or scheme $C$, let $\mathcal{O}$ denote the identity, $[m]: C \rightarrow C$ denote multiplication by $m$, $p_i: C^n \rightarrow C$ denote projection onto the $i$-th component, and $s: C^n \rightarrow C$ denote sum of all components.  For $\v \in \Z^n$, define the expression 
\begin{gather*}
D_{C,\v} = ( [v_1]\times\ldots\times[v_n])^*s^*  (\mathcal{O}) - \sum_{1 \leq k < j \leq n} v_kv_j (p_k^*\times p_j^*)s^* (\mathcal{O}) \\
- \sum_{k=1}^n \left(2v_k^2-\sum_{j=1}^nv_kv_j\right) p_k^* (\mathcal{O}) \myenspace ,
\end{gather*}
which is a divisor on the $n$-fold product $C^n$.  Over the complex numbers, the functions $\Omega_\v$ have these divisors and statisfy the elliptic net recurrence \eqref{eqn: ellrec} (see Section \ref{sec: overc}).

We now collect the results of the previous sections in one statement.
\begin{theorem}
\label{thm: arbitraryfield}
Let $n \geq 1$.  There exists a unique collection of rational functions $\Psi_{\v} \in \mathcal{K}(\mathcal{E}_S^n)$ for each $\v \in \Z^n$ such that:
\begin{enumerate}
\item The map $\v \mapsto \Psi_\v$ gives an elliptic net $W: \Z^n \rightarrow \mathcal{R}_n$.
\item $\Psi_{\v} = 1$ whenever $\v = \e_i$ for some $1 \leq i \leq n$ or $\v = \e_i + \e_j$ for some $1 \leq i < j \leq n$.
\item $\operatorname{Div}(\Psi_{\v}) = D_{\mathcal{E}_S,\v}$.
\end{enumerate}
\end{theorem}

\begin{proof}
Part (1) follows from Theorems~\ref{thm: omegav} and \ref{thm: valuation}.  Part (2) follows from Proposition~\ref{prop: 2explicit} and Theorem~\ref{thm: omegav}.  Part (3) follows from Theorem~\ref{thm: valuation}.
\end{proof}

\section{Elliptic nets from elliptic curves}
\label{sec: arbitraryfields}

%\subsection{The Elliptic Net Associated to a Curve}
%\label{sec:ellnetassociated}

In light of Theorem \ref{thm: arbitraryfield}, it is now natural to define an elliptic net associated to any cubic Weierstrass curve over any field.
\begin{definition}
\label{defn: WEP}
Let $K$ be any field.  Let $a_1, a_2, a_3, a_4, a_6 \in K$.  To this we associate a map 
\[
S = \mathbb{Z}[\alpha_1,\alpha_2,\alpha_3,\alpha_4,\alpha_6] \rightarrow K, \quad \alpha_i \mapsto a_i.
\]
Let 
$$f(x, y) = y^2 + a_1 xy + a_3 y - x^3 - a_2 x^2 - a_4 x - a_6 $$
and let $C$ be a curve defined by $f(x,y)=0$.  Then we have a Cartesian diagram
\[
\xymatrix{
\mathcal{E}_S^n \ar[d] & C^n \ar[l] \ar[d] \\
\Spec(S) & \Spec(K) \ar[l]
}
\]
under which we may pullback $\Psi_\v$ to obtain $\phi_\v \in \mathcal{K}(C^n)$ (this is possible since the fibre on the right is not contained in the support of the divisor of $\Psi_\v$, by Theorem \ref{thm: arbitraryfield}).

The non-singular points of $C$ defined over $K$, denoted $C_{ns}(K)$, form a group.  We call a set of points $\{ P_1, \ldots, P_n \}$ on the non-singular part $C_{ns}$ of a cubic curve \emph{appropriate} if the following hold:
\begin{enumerate}
\item $P_i \neq 0$ for all $i$;
\item $[2]P_i \neq 0$ for all $i$;
\item $P_i \neq \pm P_j$ for any $i \neq j$; and
\item $[3]P_1 \neq 0$ whenever $n=1$.
\end{enumerate}

If we have an appropriate $n$-tuple of points $\mathbf{P} \in C_{ns}(K)^n$, then we may define a map
$$ W_{C,\mathbf{P}} : \Z^n \rightarrow K \myenspace ,$$
defined by $W_{C,\mathbf{P}}(\v) = \phi_\v(\mathbf{P})$.  By Theorem \ref{thm: arbitraryfield}, this will be an elliptic net.  This will be called \emph{the elliptic net associated to $C$ and $\mathbf{P}$.}
\end{definition}

%Note that the divisor of $\phi_{\v}$ will be the pullback of the divisor $D_{\mathcal{E}_S,\v}$, but it may not have the form $D_{C,\v}$ because points may coincide.  For example, if $C$ is not supersingular over $\Fp$, then the divisor of $\phi_p$ over $\Fp$ has degree $p$ at each of the $p$ points of $E[p]$. {\bf Is that true?!  Doesn't the pullback take into account ramification?  Isn't the singular fibre a bigger issue?}\comm{***}

We have the following additional corollary to Theorem \ref{thm: arbitraryfield}.
\begin{corollary}
\label{cor: netzeroes}
For an elliptic net $W_{C,\mathbf{P}}: \Z^n \rightarrow K$ associated to a curve $C$ and non-singular points $\mathbf{P}$, we have $W(\v) = 0$ if and only if $\v(\mathbf{P}) = \mathcal{O}$ on $C_{ns}$.
\end{corollary}
\begin{proof}This follows from the statement that $\Omega_\v(\v \cdot \z) = 0$ if and only if $\v \cdot \z \in \Lambda$ (see Section \ref{sec: overc}).\end{proof}

\begin{example}
\label{example:net}
Figure \ref{figure: ellnet66} (in Section \ref{sec: edsedn}) shows an example elliptic net associated to the elliptic curve and points
\[
E: y^2 + y = x^3 + x^2 - 2x, \quad P=(0,0), \quad Q=(1,0)
\]
Some of the smaller terms of the net $W_{E,(P,Q)}$ can be calculated using Proposition \ref{prop: 2explicit}, for example,
\[
W(0,0) = 0, \quad W(1,0) = W(0,1) = W(1,1) = 1,
\]
\[
W(2,0) = 2y_1 + a_1x_1 + a_3 = 1, \qquad W(0,2) = 2y_2 + a_1x_2 + a_3 = 1,
\]
\[
W(1,-1) = x_2 - x_1 = 1,
\]
\[
W(2,1) = 2x_1 + x_2 - \left(\frac{y_2-y_1}{x_2-x_1}\right)^2 - a_1 \left(\frac{y_2-y_1}{x_2-x_1}\right) + a_2 = 2,
\]
\[
W(2,-1) =  (y_1+y_2)^2 - (2x_1+x_2)(x_1-x_2)^2 = -1.
\]
More terms can be calculated using the recurrence relation \eqref{eqn: ellrec}.  Since $P$ and $Q$ are independent non-torsion points, there are no zeroes in the array except the zero located at the origin ($W(0,0)=0$).  The centre row is the elliptic divisibility sequence associated to $E$ and $P$, which begins
\begin{align*}
&1, 1, -3, 11, 38, 249, -2357, 8767, 496035, -3769372, -299154043, \\& -12064147359, 632926474117,  -65604679199921, \ldots 
%\\& -6662962874355342, -720710377683595651, 285131375126739646739,\\&  5206174703484724719135, -36042157766246923788837209, 14146372186375322613610002376, \ldots
\end{align*}
The centre column is the elliptic divisibility sequence associated to $Q$.  
\end{example}

\section{Elliptic curves from elliptic nets}
\label{sec: curvesfromnets}

We are now in a position to use the results of Section \ref{sec: induction} to determine exactly which elliptic curves (or more generally cubic Weierstrass curves) give rise to any given elliptic net.

\subsection{Scale equivalence and normalisation}
\label{sec:scalenormal}

\begin{proposition}
Let $W: A \rightarrow K$ be an elliptic net.  Let $f: A \rightarrow K^*$ be a quadratic form.  Define $W^f: A \rightarrow K$ by
$$ W^f(\v) = f(\v)W(\v) \myenspace .$$
Then $W^f$ is an elliptic net.
\end{proposition}

\begin{proof}  Let $p, q, r, s \in A$.  We use multiplicative notation in $K^*$, so that $f$ satisfies
\begin{equation}
\label{eqn: quad1}
f(p+q+s)f(p)f(q)f(s)f(p+q)^{-1}f(q+s)^{-1}f(p+s)^{-1}=1 \myenspace .
\end{equation}
The parallelogram law for quadratic forms (written multiplicatively) states that
\begin{equation}
\label{eqn: quad2}
f(p-q)f(p+q) = f(p)^2f(q)^2 \myenspace .
\end{equation}
Multiplying $f(r)f(r+s)$ and equations \eqref{eqn: quad1} and \eqref{eqn: quad2} together,
\begin{gather*}
 f(p+q+s)f(p-q)f(r+s)f(r) =\\ f(q+s)f(p+s)f(r+s)f(p)f(q)f(r)f(s)^{-1} \myenspace ,
\end{gather*}
which is symmetric in $p$, $q$, and $r$, so
\begin{align*}
f(p+q+s)f(p-q)f(r+s)f(r) &= f(q+r+s)f(q-r)f(p+s)f(p) \\ 
&= f(r+p+s)f(r-p)f(q+s)f(q) \myenspace ,
\end{align*}
which shows that the recurrence \eqref{eqn: ellrec} holds for $W^f$ if it does for $W$.
\end{proof}

If two elliptic nets are related in the manner of $W$ and $W^f$ for some quadratic form $f$, then we call them \emph{scale equivalent}.  This is clearly an equivalence relation.  
\par
Let $W: \Z^n \rightarrow K$ be an elliptic net.  We say that $W$ is \emph{normalised} if $W(\e_i) = 1$ for all $1 \leq i \leq n$ and $W(\e_i + \e_j) = 1$ for all $1 \leq i < j \leq n$.  An elliptic net arising from a curve and points is normalised.  It should be emphasized that the concept of \emph{normalised} is only defined for elliptic nets with a preferred basis.

If any term of the form $W(\e_i)$, $W(2\e_i)$, $W(\e_i+\e_j)$, or $W(\e_i - \e_j)$ is zero (where $i \neq j$), or if $n=1$ and any term of the form $W(3\e_1)$ is zero, then we say that $W$ is \emph{degenerate}.  Compare the definition of \emph{degenerate} to the definition of \emph{appropriate} in Section \ref{sec: arbitraryfields}.

\begin{proposition}
\label{prop:uniqnormal}
Let $W: \Z^n \rightarrow K$ be a non-degenerate elliptic net.  Then there is exactly one scaling $W^f$ which is normalised.
\end{proposition}

\begin{proof}Define
%Specify $A_{ij} \in K^*$ for $1 \leq i \leq j \leq n$ as follows:
\begin{align*}
A_{ii} &= W(\e_i)^{-1}, \quad \mbox{ for }1 \leq i \leq n ,\\
A_{ij} &= \frac{ W(\e_i) W(\e_j)}{W(\e_i + \e_j)}, \quad \mbox{ for } 1 \le i < j \le n,\\
f(\v) &= \prod_{1 \leq i \leq j \leq n} A_{ij}^{v_iv_j} \myenspace .
\end{align*}
Then $W^f$ is normalised.  Uniqueness follows from the elementary properties of quadratic forms (as in the proof of Lemma \ref{lem:quad}).
\end{proof}

The proof demonstrates that scale equivalence has ${n+1}\choose{2}$ degrees of freedom.  If $W: \Z^n \rightarrow K$ is an elliptic net, then its \emph{normalisation} $\widetilde W$ is defined to be the unique normalised elliptic net which is a scaling of $W$.  A \emph{coordinate sublattice} of $\Z^n$ is a sublattice of the form 
$$\{ \v \in \Z^n : v_i = 0 \mbox{ for }i \notin I \}$$ 
for some proper non-empty subset $I \subset \{1,2,\ldots,n\}$.  The \emph{rank} of the sublattice is $\#I$.

\subsection{Curves from nets of ranks $1$ and $2$}

Define a change of variables of a cubic curve in Weierstrass form to be \emph{unihomothetic} if it is of the form
\begin{equation}
\label{eqn: coordchange}
x' = x + r , \quad y' = y + sx + t,
\end{equation}
for some $r$, $s$ and $t$.

The rank-one result in the following form is due to Christine Swart.

\begin{proposition}[{Swart \cite[Theorem 4.5.3]{Swa}}]
\label{prop: curvefromnet1}
Let $W: \Z \rightarrow K$ be a normalised non-degenerate elliptic net.  Then the family of curve-point pairs $(C,P)$ such that $W = W_{C,P}$ is three dimensional.  These are the curve and non-singular point
$$C: y^2 + a_1xy + a_3y = x^3 + a_2x^2 + a_4x + a_6, \quad P = (0,0), $$
where
\begin{align*}
a_1 &= \frac{W(4) + W(2)^5 - 2W(2)W(3)}{W(2)^2W(3)}, \\
a_2 &= \frac{W(2)W(3)^2 + (W(4) + W(2)^5) - W(2)W(3)}{W(2)^3W(3)},
\end{align*}
\[
a_3 = W(2), \qquad a_4 = 1, \qquad a_6 = 0,
\]
or any image of these under a unihomothetic change of coordinates.
\end{proposition}

\begin{proof}See \cite[Theorem 4.5.3]{Swa}.  A normalised rank 1 non-degenerate elliptic net has $W(2)~\neq~0$ and $W(3)~\neq~0$.  Any singular point $P=(x,y)$ on a cubic Weierstrass curve has vanishing partial derivatives, which implies $\Psi_2(P)~=~2y~+~a_1x~+~a_3~=~0$ (see Proposition \ref{prop: 2explicit}).  Therefore, if any curve and singular point gives rise to $W$, then $W(2)~=~0$, in contradiction to non-degeneracy.  The division polynomials $\Psi_1$, $\Psi_2$, $\Psi_3$ and $\Psi_4$ are invariant under a change of coordinates of the form \eqref{eqn: coordchange}.  Then, it is a simple calculation to check that $W_{C,P}$ agrees with $W$ at the first four terms; hence $W_{C,P} = W$ by Theorem \ref{thm: 1laurent}.  Conversely, suppose $W = W_{C', P'}$.  After applying a transformation of the form \eqref{eqn: coordchange} taking $P'$ to $(0,0)$ and taking $a_4$ to $1$, substitution of the division polynomials into the equations above verifies that $a_i' = a_i$ for all $i$. 
\end{proof}

\begin{proposition}
\label{prop: curvefromnet2}
Let $W: \Z^2 \rightarrow K$ be a normalised non-degenerate elliptic net.  Then the family of $3$-tuples $(C,P_1,P_2)$ such that $W = W_{C,P_1,P_2}$ is three dimensional.  These are the curve and non-singular points
$$C: y^2 + a_1xy + a_3y = x^3 + a_2x^2 + a_4x + a_6, $$
$$ P_1 = (0,0), \quad P_2 = (W(1,2) - W(2,1),0) \myenspace ,$$
with
$$ a_1 = \frac{W(2,0) - W(0,2)}{W(2,1)-W(1,2)}, \quad a_2 = 2W(2,1)- W(1,2), \quad a_3 = W(2,0)$$
$$a_4 = (W(2,1) - W(1,2)) W(2,1), \qquad a_6 = 0, $$
or any image of these under a unihomothetic change of coordinates.
\end{proposition}

\begin{proof}
In a normalised non-degenerate elliptic net, 
\[
W(2,1)-W(1,2)=W(1,-1)\neq0, \quad W(2,0) \neq 0, \quad W(0,2) \neq 0
\]
(see Theorem \ref{thm: 2laurent}).  Thus (as in the previous theorem) if a curve and points give rise to $W$, then the points are non-singular.  The formul{\ae} for $W(2,0)$, $W(0,2)$, $W(2,1)$ and $W(1,2)$ are invariant under a change of coordinates of the form \eqref{eqn: coordchange}.  The net $W_{C,P_1,P_2}$ agrees with $W$ at the terms $(2,0)$, $(0,2)$, $(2,1)$ and $(1,2)$; hence $W_{C,P_1,P_2} = W$ by Theorem \ref{thm: 2laurent}.  Conversely, suppose $W = W_{C', P_1',P_2'}$.  After applying a unihomothetic transformation taking $P_1'$ to $(0,0)$ and $P_2'$ to $(W(1,2) - W(2,1),0)$, substitution of the net polynomials into the equations above verifies that $a_i' = a_i$ for all $i$.
\end{proof}

\begin{example}
Plugging in terms from the elliptic net of Figure \ref{figure: ellnet66} (in Section \ref{sec: edsedn}) to the formul{\ae} in the statement of Proposition \ref{prop: curvefromnet2} we recover the curve and points stated in the figure caption.
\end{example}

\begin{remark}
A more symmetric set of equations in the case of characteristic not equal to $2$ is as follows:
\[
P_1 = (v,0) , \quad P_2 = (-v,0), \quad 2v = W(2,1) - W(1,2) 
 \myenspace ,
\]
\[
a_1 = \frac{W(2,0) - W(0,2)}{W(2,1)-W(1,2)}, \quad
2a_2 = W(2,1) + W(1,2),
\]
\[
2a_3 = W(2,0) + W(0,2), \quad 4a_4 = -(W(2,1) - W(1,2))^2,
\]
\[
8a_6 = - (W(2,1) - W(1,2))^2(W(2,1) + W(1,2)).
\]
\end{remark}

\subsection{Curves from nets in general rank}

\begin{theorem}
\label{thm: curvefromnetn}
Let $n \ge 1$.  Let $W: \Z^n \rightarrow K$ be a normalised non-degenerate elliptic net.  Then the set of curves $C$ and $\mathbf{P} \in C^n$ such that $W = W_{C,\mathbf{P}}$ forms a three-dimensional family of tuples $(C,\mathbf{P})$.  Further, none of the points $P \in \mathbf{P}$ are singular.  In particular, the family consists of one such tuple and all its images under unihomothetic changes of coordinates.
\end{theorem}

\begin{proof}
The proof is by strong induction on $n$, where the inductive statement has two parts:  (I) that the theorem holds for $n$; and (II) that $W(\v) \neq 0$ for some $\v \in \{\pm 1 \}^n$.  The base case consists of ranks $n=1,2$:  part (I) is by Propositions \ref{prop: curvefromnet1} and \ref{prop: curvefromnet2}; part (II) is by non-degeneracy, which implies $W(\e_1) \neq 0$ and $W(\e_1 + \e_2) \neq 0$.

Suppose $n \geq 3$ and the inductive statement holds for all $k < n$.  Let $W_1$, $W_2, \ldots, W_n$ be the normalised elliptic subnets of $W$ associated to the rank $n-1$ coordinate sublattices $L_i = \{ \mathbf{v} : v_i = 0 \}$.  These are defined as nets $W_i: L_i \rightarrow K$ but they can be identified with nets $W_i': \Z^{n-1} \rightarrow K$ in the obvious way (by deleting the zero coordinate).  They are normalised and non-degenerate (by definition, non-degeneracy at rank $n$ implies non-degeneracy on rank $n-1$ sublattices for $n > 2$).  By the inductive hypothesis part (I), we have
$W_i' = W_{C_i, \mathbf{P}_i}$
for some curves $C_i$ and non-singular points $\mathbf{P}_i \in C_i^{n-1}$.

We observe a consequence of Proposition \ref{prop: psiuuomega}.  Suppose $V_1:\Z^m \rightarrow K$ is an elliptic net of rank $m$ associated to $C$ and $\mathbf{P}$.  Also suppose $$V_2:\{ \v \in \Z^m : v_m = 0 \} \rightarrow K$$ is the elliptic subnet of $V_1$ associated to the coordinate sublattice of rank $m-1$ which consists of vectors with last coordinate zero.  Suppose $V_2':\Z^{m-1} \rightarrow K$ is naturally identified with $V_2$ by simply deleting the last coordinate of the domain.  Then $V_2'$ is associated to $C$ and $\mathbf{P}'$ where $\mathbf{P}'$ is simply $\mathbf{P}$ with the last coordinate deleted.  This statement, appropriately adjusted, holds for any coordinate hyperplane (not just the one with last coordinate zero).

Consider two of the rank $n-1$ subnets, say $W_i$ and $W_j$.  Let $W_{ij} = W_i \cap W_j$ in $W$.  Define $W_{ij}' : \Z^{n-2} \rightarrow K$ by the obvious identification.  Then, $W_{ij}' = W_{C_{ij},\mathbf{P}_{ij}}$ for some curve $C_{ij}$ and $\mathbf{P}_{ij} \in C_{ij}^{n-2}$. By the foregoing, $C_i = C_j = C_{ij}$, $\mathbf{P}_{ij}$ is just $\mathbf{P}_j$ with the $i$-th coordinate deleted, and $\mathbf{P}_{ij}$ is just $\mathbf{P}_i$ with the $(j-1)$-th coordinate deleted.

Considering every such pair, we may define a candidate curve $C$ by $C = C_i$ for all $i$ and $\mathbf{P} \in C^n$ defined as the unique $n$-tuple which results in $\mathbf{P}_i$ upon deleting the $i$-th coordinate.  By the foregoing, this is well-defined.  Now we see that $W$ agrees with $W_{C,\mathbf{P}}$ on all coordinate sublattices of rank $n-1$.  By the inductive statement part (II) and Theorem \ref{thm: nlaurent}, we see that $W$ is determined by its sublattices of rank $n-1$.  Therefore $W = W_{C,\mathbf{P}}$.  

To show part (II) of the inductive statement, we observe that if $W(\v)=0$ for all $\v \in \{ \pm 1 \}^n$, then $\v(\mathbf{P}) = \mathcal{O}$ for all such $\v$ (by Corollary \ref{cor: netzeroes}).  But this is impossible, since it would imply $[2]P_i = \mathcal{O}$ for $1 \le i \le n$, a contradiction to non-degeneracy (again Corollary \ref{cor: netzeroes}).

A change of coordinates of the form \eqref{eqn: coordchange} for $C$ does not change the elliptic net, as it is determined by its values on its coordinate hyperplanes, where this is true.  Further, if two tuples \emph{not} related by such a change of coordinates generate the same net $W$, then the same would hold for some coordinate hyperplane, a contradiction.  This demonstrates part (I) of the inductive statement.\end{proof}

\section{The curve-net theorem}
\label{sec: curvenet}

We set some remaining terminology, and then proceed to the statement and proof of the main theorem.

\subsection{Homothety and singular elliptic nets}
\label{sec:homosing}

The only changes of coordinates of a Weierstrass equation into another are compositions of unihomothetic changes of coordinates and changes of coordinates of the form $(x,y) \mapsto (\lambda^2x, \lambda^3y)$, which we refer to as \emph{homotheties} (since they correspond to homotheties of the lattice in the complex uniformization).
\begin{proposition}
\label{prop: homothety}
Consider the rank $n$ elliptic net $W_{C,\mathbf{P}}$ associated to 
$$C: y^2 + a_1 xy + a_3 y= x^3 + a_2 x^2 + a_4 x + a_6$$
defined over $K$ and $\mathbf{P} \in C(K)^n$.  Let $\lambda$ be a non-zero element of $K$.  Suppose $\phi_\lambda: C \rightarrow C_\lambda$ is the isomorphism of curves taking $C$ to
$$C_\lambda: y^2 + \lambda a_1 xy + \lambda^3 a_3 y = x^3 + \lambda ^2 a_2 x^2 + \lambda^4 a_4 x + \lambda^6 a_6$$
under the change of coordinates $(x,y) \mapsto (\lambda^2x, \lambda^3y)$.  Then
$$\widetilde W_{C_\lambda, \phi_\lambda(\mathbf{P})} = \lambda \widetilde W_{C,\mathbf{P}}$$
In particular, let $\delta_{ij}$ be the Kronecker delta, and define
$$ g(\v) = {-1 - \sum_{1 \leq i < j \leq n} (-1)^{\delta_{ij}} v_iv_j} \myenspace .$$
Then
$$W_{C_\lambda, \phi_\lambda(\mathbf{P})} = \lambda^{g(\v)}W_{C, \mathbf{P}} \myenspace .$$
\end{proposition}

\begin{proof}
The first statement is entailed by the second.  From the general theory of Weierstrass sigma functions, $\sigma( \lambda z, \lambda \Lambda) = \lambda \sigma(z, \Lambda)$.  Thus by Definition \ref{defn: psi}, we know that $\Omega_\v(\lambda \z; \lambda \Lambda) = \lambda^{g(\v)}  \Omega_\v(\z; \Lambda)$.  As in Section \ref{sec: netpoly}, this allows us to conclude that the same holds for $\Psi_\v$, so that
\[
\Psi_\v(\lambda^2x, \lambda^3 y, \lambda^i \alpha_i) = \lambda^{g(\v)}  \Psi_\v(x, y, \alpha_i),
\]
from which the result follows. 
%If $\lambda$ is a non-zero element of $K$ then $C_\lambda$ is again an elliptic curve.  
%The proposition holds for $\Omega_\v$ over $\C$ by Definition \ref{defn: psi}.  Therefore the rational function representations of Theorem \ref{thm: arbitraryfield} are weighted homogeneous in the appropriate way; hence it holds over any field.  Is this any kind of proof????****
\end{proof}

Therefore we set the following definition

\begin{definition}
If $W : \Z^n \rightarrow K$ is an elliptic net, then with the notation of Proposition \ref{prop: homothety}, we define
$$W^\lambda (\v) \mathcolon= \lambda^{g(\v)} W(\v) \myenspace .$$
This gives an action of $K$ on elliptic nets $W: \Z^n \rightarrow K$ called the \emph{homothety action}.  Two elliptic nets are \emph{homothetic} if they are in the same orbit of the action of $K$.
\end{definition}

The following proposition is immediate.

\begin{proposition}
Let $W: \Z^n \rightarrow K$ be an elliptic net.  Then for any non-zero $\lambda \in K$, $W^\lambda$ is normalised if and only if $W$ is.
\end{proposition}

Let $W: \Z^n \rightarrow K$ be an elliptic net.  If the curve $C$ associated to its normalisation is a nodal or cuspidal cubic, then $W$ is called \emph{singular}.  If, instead, $C$ is an elliptic curve, then $W$ is called \emph{non-singular}.  In either case, the discriminant $\Delta$ of $W$ is defined to be the discriminant of the associated Weierstrass equation.  Similarly, the $j$-invariant is the $j$-invariant of the associated Weierstrass equation.  The discriminant of an elliptic net changes by a factor of $\lambda^{12}$ under homothety, while the $j$-invariant remains unaltered.

\subsection{The curve-net theorem}

We may put a partial ordering on tuples $(C,P_1,\ldots,P_n)$ where $C$ is a Weierstrass curve and $P_i$ are non-singular points on the curve.  We do this by defining 
\[
(C, P_1, \ldots, P_n) \leq (D, Q_1, \ldots, Q_m)
\]
if and only if $C=D$ and the groups they generate satisfy a containment 
\[
\left<P_1, \ldots, P_n \right> \subseteq \left<Q_1,\ldots,Q_n \right>.
\]
The collection of all elliptic nets is partially ordered by the subnet relation.  Collecting our work up to this point, we have now shown:
\begin{theorem}
\label{thm:big}
For each field $K$, there is an explicit isomorphism of partially ordered sets
\begin{equation*}
\xymatrix{
\left\{
{
\begin{array}{l}
\mbox{scale equivalence classes of } \\
\mbox{non-degenerate elliptic nets}\\
W: \Z^n \rightarrow K
\mbox{ for some }n
\end{array}
}
\right\} \ar[d] \\
\left\{
{
\begin{array}{l}
\mbox{tuples }(C, P_1, \ldots , P_m )\mbox{ for some }m, \mbox{ where $C$}\\
\mbox{is a cubic curve in Weierstrass form over $K$,}\\
\mbox{considered modulo unihomothetic changes}\\
\mbox{of variables, and such that }\{P_i\} \in C_{ns}(K)^m \\
\mbox{is appropriate}
\end{array}
}
\right\} \ar[u]
}
\end{equation*}
Non-singular nets correspond to elliptic curves.  The action of $K$ (by homothety) on the sets preserves the order and respects the isomorphism.  The bijection takes an elliptic net of rank $n$ to a tuple with $n$ points.  The elliptic net $W$ associated to a tuple $(C,P_1,\ldots,P_n)$ satisfies the property that $W(v_1,\ldots,v_n) = 0 $ if and only if $v_1P_1 + \ldots + v_nP_n = 0$ on the curve $C$.
\end{theorem}

\begin{proof} In the diagram in the statement of the theorem, call the upper set $\mathcal{N}$ and the lower set $\mathcal{C}$.  The first claim is that there is an injective map $\mathcal{N} \rightarrow \mathcal{C}$.  Proposition \ref{prop:uniqnormal} shows that each scale equivalence classes in $\mathcal{N}$ contains a unique normalised elliptic net, so we can define the map by Theorem \ref{thm: curvefromnetn} (which also guarantees injectivity).  Corollary \ref{cor: netzeroes} shows that the result is an element of $\mathcal{C}$.  This shows the first claim.
\par
The second claim is that there exists an inverse map $\mathcal{C} \rightarrow \mathcal{N}$.  The map is given by Definition \ref{defn: WEP}, which is well-defined as a result of Theorem \ref{thm: arbitraryfield}.  It is required to check that the resulting elliptic net is normalised (Proposition \ref{prop: 2explicit}) and non-degenerate (Corollary \ref{cor: netzeroes}).  Theorem \ref{thm: curvefromnetn} says that this map is indeed an inverse to the map of the first claim.  This gives the second claim and the bijection of sets.
\par
It is clear that the bijection associates an elliptic net of rank $n$ to a tuple with $n$ points, and that it preserves the partial ordering.  The action of homothety is preserved by Proposition \ref{prop: homothety}.  And the final statement of the theorem is a result of Corollary \ref{cor: netzeroes}.
\end{proof}

\begin{remark}
The degenerate cases present several difficulties.  One is that a degenerate elliptic net may not be determined by the usual initial set of terms as given in Section \ref{sec: induction}.  For example, the sequence
\[
W(n) = \left\{ \begin{array}{ll}
0 & n \neq k\\
1 & n = k
\end{array} \right. .
\]
is an elliptic net for any non-zero integer $k$.  However, some degenerate sequences can be thought of as arising from singular points on a singular cubic.  For example, consider a sequence associated to an elliptic curve $E$ and point $P$ both defined over $\Q$ such that $P$ reduces to a singular point modulo some prime $p$.  Then the sequence regarded modulo $p$ as living in $\Fp$ (which is necessarily a degenerate elliptic net) should be associated to a point on the special fibre of the N\'eron model.  It is likely that Theorem \ref{thm:big} can be extended to include these cases (this is future work).
\end{remark}

\bibliographystyle{plain}

\begin{thebibliography}{10}

\bibitem{Aya}
Mohamed Ayad.
\newblock P\'eriodicit\'e (mod {$q$}) des suites elliptiques et points
  {$S$}-entiers sur les courbes elliptiques.
\newblock {\em Ann. Inst. Fourier (Grenoble)}, 43(3):585--618, 1993.

\bibitem{Cha}
K.~Chandrasekharan.
\newblock {\em Elliptic functions}, volume 281 of {\em Grundlehren der
  Mathematischen Wissenschaften [Fundamental Principles of Mathematical
  Sciences]}.
\newblock Springer-Verlag, Berlin, 1985.

\bibitem{ChuChu}
D.~V. Chudnovsky and G.~V. Chudnovsky.
\newblock Sequences of numbers generated by addition in formal groups and new
  primality and factorization tests.
\newblock {\em Adv. in Appl. Math.}, 7(4):385--434, 1986.

\bibitem{CorZah}
Gunther Cornelissen and Karim Zahidi.
\newblock Elliptic divisibility sequences and undecidable problems about
  rational points.
\newblock {\em J. Reine Angew. Math.}, 613:1--33, 2007.

\bibitem{EisEve}
Kirsten Eisentr{\"a}ger and Graham Everest.
\newblock Descent on elliptic curves and {H}ilbert's tenth problem.
\newblock {\em Proc. Amer. Math. Soc.}, 137(6):1951--1959, 2009.

\bibitem{EveMclWar}
Graham Everest, Gerard Mclaren, and Thomas Ward.
\newblock Primitive divisors of elliptic divisibility sequences.
\newblock {\em J. Number Theory}, 118(1):71--89, 2006.

\bibitem{EveMilSte}
Graham Everest, Victor Miller, and Nelson Stephens.
\newblock Primes generated by elliptic curves.
\newblock {\em Proc. Amer. Math. Soc.}, 132(4):955--963 (electronic), 2004.

\bibitem{EvePooShpWar}
Graham Everest, Alf van~der Poorten, Igor Shparlinski, and Thomas Ward.
\newblock {\em Recurrence Sequences}, chapter Elliptic Divisibility Sequences,
  pages 163--175.
\newblock American Mathematical Society, Providence, 2003.

\bibitem{FomZel}
Sergey Fomin and Andrei Zelevinsky.
\newblock The {L}aurent phenomenon.
\newblock {\em Adv. Appl. Math.}, 28(2):119--144, 2002.

\bibitem{FreLan}
Gerhard Frey and Tanja Lange.
\newblock Background on curves and {J}acobians.
\newblock In {\em Handbook of elliptic and hyperelliptic curve cryptography},
  Discrete Math. Appl. (Boca Raton), pages 45--85. Chapman \& Hall/CRC, Boca
  Raton, FL, 2006.

\bibitem{GasRah}
George Gasper and Mizan Rahman.
\newblock {\em Basic hypergeometric series}, volume~96 of {\em Encyclopedia of
  Mathematics and its Applications}.
\newblock Cambridge University Press, Cambridge, second edition, 2004.
\newblock With a foreword by Richard Askey.

\bibitem{Hon}
A.~N.~W. Hone.
\newblock Elliptic curves and quadratic recurrence sequences.
\newblock {\em Bull. London Math. Soc.}, 37(2):161--171, 2005.

\bibitem{Ing2}
Patrick Ingram.
\newblock Multiples of integral points on elliptic curves.
\newblock {\em J. Number Theory}, 129(1):182--208, 2009.

\bibitem{MazTat}
B.~Mazur and J.~Tate.
\newblock The {$p$}-adic sigma function.
\newblock {\em Duke Math. J.}, 62(3):663--688, 1991.

\bibitem{Poo2}
Bjorn Poonen.
\newblock Hilbert's tenth problem and {M}azur's conjecture for large subrings
  of {$\Bbb Q$}.
\newblock {\em J. Amer. Math. Soc.}, 16(4):981--990 (electronic), 2003.

\bibitem{Rob}
James Propp.
\newblock Robbins forum.
\newblock \url{http://jamespropp.org/about-robbins}.

\bibitem{Shi}
Rachel Shipsey.
\newblock {\em Elliptic Divibility Sequences}.
\newblock PhD thesis, Goldsmiths, University of London, 2001.

\bibitem{Sil5}
Joseph~H. Silverman.
\newblock Common divisors of elliptic divisibility sequences over function
  fields.
\newblock {\em Manuscripta Math.}, 114(4):431--446, 2004.

\bibitem{Sil4}
Joseph~H. Silverman.
\newblock {$p$}-adic properties of division polynomials and elliptic
  divisibility sequences.
\newblock {\em Math. Ann.}, 332(2):443--471 (Addendum 473--474), 2005.

\bibitem{Sil1}
Joseph~H. Silverman.
\newblock {\em The arithmetic of elliptic curves}, volume 106 of {\em Graduate
  Texts in Mathematics}.
\newblock Springer, Dordrecht, second edition, 2009.

\bibitem{Sta3}
Katherine~E. Stange.
\newblock The {T}ate pairing via elliptic nets.
\newblock In {\em Pairing-Based Cryptography - PAIRING 2007}, volume 4575 of
  {\em Lecture Notes in Comput. Sci.}, pages 329--348. Springer, Berlin, 2007.

\bibitem{Str}
Marco Streng.
\newblock Divisibility sequences for elliptic curves with complex
  multiplication.
\newblock {\em Algebra Number Theory}, 2(2):183--208, 2008.

\bibitem{Swa}
Christine Swart.
\newblock {\em Elliptic curves and related sequences}.
\newblock PhD thesis, Royal Holloway and Bedford New College, University of
  London, 2003.

\bibitem{Poo}
Alfred~J. van~der Poorten.
\newblock Elliptic curves and continued fractions.
\newblock {\em J. Integer Seq.}, 8(2):Article 05.2.5, 19 pp. (electronic),
  2005.

\bibitem{PooSwa}
Alfred~J. van~der Poorten and Christine~S. Swart.
\newblock Recurrence relations for elliptic sequences: every {S}omos 4 is a
  {S}omos {$k$}.
\newblock {\em Bull. London Math. Soc.}, 38(4):546--554, 2006.

\bibitem{War}
Morgan Ward.
\newblock Memoir on elliptic divisibility sequences.
\newblock {\em Amer. J. Math.}, 70:31--74, 1948.

\bibitem{WenEkhCha}
Chu Wenchang, Shalosh~B. Ekhad, and Robin~J. Chapman.
\newblock Problems and {S}olutions: {S}olutions: 10226.
\newblock {\em Amer. Math. Monthly}, 103(2):175--177, 1996.

\end{thebibliography}
\def\cprime{$'$}

\end{document}